%% file: Dissertation.tex
\documentclass[economy]{dukedissertation}

\usepackage{amsmath, amssymb, amsfonts, amsthm, dsfont}
\usepackage{graphicx}
\usepackage{color}
\usepackage{bm}
\usepackage{subfigure}
\usepackage{graphicx}
\usepackage{mathabx}
\usepackage{multirow}
\usepackage{setspace}
							
\newtheorem{theorem}{Theorem}
\newtheorem{lemma}[theorem]{Lemma}

\def\dsst{\displaystyle}
\def\R{\mathds{R}}
\theoremstyle{definition}
\newtheorem{definition}{Definition}[section]

\author{Ching-Yin Wong}
\title{Sigma Models with Repulsive Potentials}
\department{Mathematics} 

\date{2017} 

\copyrighttext{ All rights reserved except the rights granted by the\\
   \href{http://creativecommons.org/licenses/by-nc/3.0/us/}
        {Creative Commons Attribution-Noncommercial Licence}
}

\usepackage[pdftex, pdfusetitle, plainpages=false, 
				letterpaper, bookmarks = false, bookmarksnumbered,
				colorlinks, linkcolor=black, citecolor=black,
	         filecolor=black, urlcolor=black]
				{hyperref}

\begin{document}

\maketitle

\include{{./Abstract/abstract}}

\dedication{To my parents}

\tableofcontents 
\include{{./Abbreviations/listofabbr}} 


%
%
%
\include{{./Chapter1/introduction}}
\include{{./Chapter2/Chapter2}} 
\include{{./Chapter3/Chapter3}}
\include{{./Chapter4/Chapter4}}
\include{{./Chapter5/Chapter5}}



\nocite{*} 
\cleardoublepage
\normalbaselines 
\bibliography{./Bibliography/mybib} 


\end{document}

%% file: Abstract/abstract.tex
\abstract

Motivated by questions arising in the study of harmonic maps and Yang Mills theory, we study new techniques for producing optimal monotonicity relations for geometric partial differential equations. We apply these results to sharpen epsilon regularity results. As a sample application, we analyze energy minimizing maps from compact manifolds to the space of Hermitian matrices, where the energy of the map includes the usual kinetic term and a singular potential designed to force the image of the map to lie in a set homotopic to a Grassmannian.

%% file: Abbreviations/listofabbr.tex
\abbreviations


\section*{Symbols}


\begin{symbollist}
	\item[$\mathds{R}^k$] $k$-dimensional Euclidean space
	\item[$M, N$] Riemannian manifolds
	\item[$\dim(M)$] Dimension of $M$
	\item[$g_{ij},\ \tilde{g}_{\alpha\beta}$] Metric on Riemannian manifolds
	\item[$\Gamma_{ij}^k,\ \tilde{\Gamma}_{\alpha\beta}^\gamma$] Levi-Civita connection on Riemannian manifolds
	\item[$TM$] Tangent bundle of $M$
	\item[$T^*M$] Cotangent bundle of $M$ 
	\item[$\textrm{Ric}$] Ricci curvature 2-tensor
	\item[$\textrm{Hess}(W)$] Hessian matrix of some function $W$
	\item[$d,\ d^*$] Exterior derivative and its adjoint
	\item[$d\textrm{vol}$] Volume form of some Riemannian manifold
	\item[$\Delta$] Laplace-Beltrami operator
	\item[$\nabla_X$] Covariant derivative along vector field $X$
	\item[$B_R(x)$] Geodesic ball centered at $x$ with radius $R$
	\item[$S_R(x)$] Geodesic sphere centered at $x$ with radius $R$
	\item[$P_R(x,t)$] The set $\left\{(y,\tau)\in M| |y-x|\leq R,\ |\tau-t|\leq R^2 \right\}$ of Riemannian manifold $M$.
	\item[$T_R(t)$] The set $M\times [t-4R^2, t-R^2]$ of Riemannian manifold $M$
\end{symbollist}




%% file: Chapter1/introduction.tex
\chapter{Introduction}

\section{Harmonic Maps}

Let $M$ be a compact Riemannian manifold and $N$ a complete Riemannian manifold. Let $f:M\rightarrow N$ be a smooth map. $df$ is a section of the vector bundle $T^*M\otimes f^*(TN)$ over $M$. Denote the metrics of $M$ and $N$ by  $g^{ij}$ and $ \tilde{g}^{\alpha\beta}$ and the respective connections by  $\Gamma_{ij}^k$ and $ \widetilde{\Gamma}_{\beta\gamma}^\alpha$. The kinetic energy functional is defined to be
\begin{equation} 
E(f)=\int_M \frac{1}{2}|df|^2\, d\text{vol}=\int_M e(f)\, d\text{vol}. \label{energy}
\end{equation}

We are interested in finding representatives of each homotopy class of $M$ that minimize the above energy. Harmonic maps are the critical points of $E(f)$. They satisfy the corresponding Euler-Lagrange equation:
\begin{equation}\label{ellipticharmonic} 
\Delta f^\alpha+g^{ij}\widetilde{\Gamma}_{\beta\gamma}^\alpha\frac{\partial f^\beta}{\partial x^i}\frac{\partial f^\gamma}{\partial x^j}=0,\qquad \quad \forall \alpha=1,\cdots,\dim(N). 
\end{equation}
Here, we follow the Einstein summation notation. $\Delta $ is the Laplace-Beltrami operator on $M$. When the dimension of $M$ is 1, harmonic maps are closed geodesics on $N$. \\

This Euler-Lagrange equation is a non-linear elliptic equation, where the nonlinearity arises from the quadratic term
\begin{equation}
Q^\alpha(df, df)=g^{ij}\widetilde{\Gamma}_{\beta\gamma}^\alpha\frac{\partial f^\beta}{\partial x^i}\frac{\partial f^\gamma}{\partial x^j}.
\end{equation}

A typical approach to constructing harmonic maps is to study the solution $f_t$ to the corresponding parabolic flow
\begin{equation}\label{parabolicharmonic}
\frac{\partial f}{\partial t}+\Delta f+ Q(df,df)=0.
\end{equation} 

This is the gradient flow of $E(f)$, which decreases $E(f)$. We want to know whether $f_t$ converges to a harmonic map when $t\rightarrow \infty$. One of the difficulties comes from the quartic term of the induced equation of the energy density 
\begin{equation}
e(f) = \frac{1}{2}|df|^2.
\end{equation}

The sign of this quartic term is determined by the Riemannian curvature of the target manifold $N$.\\

The following result from Eells and Sampson \cite{EellsSampson} shows the existence of harmonic maps in each homotopy class of $M$ when the target manifold $N$ has non-positive Riemannian curvature:
\begin{theorem}[Eells - Sampson \cite{EellsSampson}]
	Let $M$ be compact and $N$ complete. Suppose $N$ has non-positive Riemannian curvature and $f:M\rightarrow N$ is a continuous map. Let $f_t$ be the solution to \eqref{parabolicharmonic} with $f_t(0,x)=f_0(x)$. If $f_t$ is bounded as $t\rightarrow\infty$, then $f_t$ is homotopic to a harmonic map $f_\infty$ with $E(f_\infty)\leq E(f)$. Furthermore, if $N$ is compact with non-positive Riemannian curvature, every homotopy class of maps from $M$ to $N$ contains a harmonic map that minimizes the energy $E(f)$. 
\end{theorem}

While negatively curved target manifolds guarantee the existence of harmonic maps, the existence of harmonic maps when $f$ is mapped into a non-negatively curved manifold is more complicated. For example, Schoen and Uhlenbeck \cite{S-USphere} proved the following result:
\begin{theorem}[Schoen-Uhlenbeck \cite{S-USphere}]
	Consider maps $f:\R^n\rightarrow S^k$, where $S^k$ is the Euclidean $k$-sphere. Let
	$$ d(k)=\left\{\begin{array}{ll} 3, & k=3 \\ \min\{\frac{1}{2}k+1, 6\}, & k\geq 4 \end{array} \right. $$
	If $k\geq 3$ and $n\leq d(k)$, there is no non-constant minimizing harmonic map $f$ from $R^n$ to $S^k$.  
\end{theorem} 

Leung (\cite{leung1982stability}) showed the following
\begin{theorem}[Leung \cite{leung1982stability}]
	For $n\geq 3$, there exists no non-constant stable harmonic maps from any compact Riemannian manifolds to $S^n$. 
\end{theorem} 

Since minimizing harmonic maps may not always exist for a given target manifold, we are interested in understanding when  smooth harmonic maps do exist and how singularities develop when smooth minimal harmonic maps do not exist. 
In the elliptic case, Schoen and Uhlenbeck \cite{S-UHarmonic} showed that for any energy minimizing map $f:M\rightarrow N$ with $\left\|f\right\|_{W^{1,2}}<\infty$ and $\text{image}(f)$ is a compact subset of $N$, the Hausdorff co-dimension of its singularity set is at least 3. 

\section{Sigma Model Approach}

Instead of considering the usual energy functional $\dsst E(f)=\int_M \frac{1}{2}|df|^2\, d\text{vol} $, which includes only the kinetic energy, we first isometrically embed $N$ in $R^k$ and now allow the image of $f$ to stay in a tubular neighborhood of the target manifold $N$, i.e., $f:M\rightarrow R^{k}\supset N$. Imposing a potential energy $W(f)$, we hope the solution will be pushed to stay on the target manifold. Therefore, we consider
\begin{equation}
E(f) = \int_M\left(\frac{1}{2}|df|^2+W(f)\right)\, d\textrm{vol}.
\end{equation}  

We denote the new energy density as
\begin{equation}
e(f) = \frac{1}{2}|df|^2+W(f).
\end{equation}

By variational calculus,
\begin{equation}
\begin{split}
\frac{d}{ds}E(f_s)\Big|_{s=0} = & \int_M \left(\langle df, d\left(\frac{\partial f}{\partial s}\right) \rangle + \langle \nabla W(f), \frac{\partial f}{\partial s}\rangle\right) \, d\textrm{vol} \Big|_{s=0}\\
= & \int_M \left(\langle d^*df+\nabla W(f), \frac{\partial f}{\partial s}\rangle\right) \, d\textrm{vol} \Big|_{s=0}.
\end{split}
\end{equation}

Thus, $f$ is a critical point of $E$ only when $f$ satisfies the elliptic equation
\begin{equation}
d^*df+\nabla W(f) = 0.
\end{equation}

The induced Laplace equation of the energy density $e(f)$ of the critical map $f$ is
\begin{equation}
\begin{split}
\Delta e(f) = & \Delta \left(\frac{1}{2}|df|^2+W(f)\right)\\
= & -|\nabla df|^2+\langle dd^*df, df\rangle-\text{Ric}(df, df)+\langle \nabla W, d^*df\rangle - \langle d\nabla W(f), df\rangle \\
= & -|\nabla df|^2 -\text{Ric}(df, df)-|\nabla W|^2-2\text{Hess}(W)_{AB}\langle df^A, df^B\rangle.
\end{split}
\end{equation}

In the last step, we apply the fact that $f$ satisfies $d^*df+\nabla W(f)=0$. $\text{Hess}(W)_{AB}$ is the Hessian matrix of $W$ in the $A,B$ direction.

We also examine the corresponding gradient flow. Now assume $f(t, x):M\times I\rightarrow R^k\supset N$ is a family of functions such that $f(t, x)$ (also denoted as $f_t(x)$) satisfies 
\begin{equation}\label{gradientflow}
\frac{\partial f_t}{\partial t} + d^*df_t+\nabla W(f_t) = 0,\qquad \text{with $f(0, x)=f_0(x)$ smooth and its image on $N$.}
\end{equation}

It can be easily shown that the energy $E(f_t)$ decreases along the flow
\begin{equation}
\frac{d}{dt}E(f_t) = \int_M \left(\langle d^*df_t+\nabla W(f_t), \frac{\partial f_t}{\partial t}\rangle\right) \, d\textrm{vol} = -\int_M \left|\frac{\partial f_t}{\partial t}\right|^2\, dx \leq 0.
\end{equation}

Our focus has turned to the existence and regularity of $f_t(x)$ when $t\rightarrow \infty$ and if so, whether $f_\infty(x) = \dsst \lim_{t\rightarrow\infty} f_t(x)$ will stay on the target manifold $N$. 

The induced heat equation of the energy density $e(f_t)$ of \eqref{gradientflow} is
\begin{equation}
\begin{split}
\left(\frac{\partial}{\partial t}+\Delta\right)e(f_t) = & \left(\frac{\partial}{\partial t}+\Delta\right)\left(\frac{1}{2}|df_t|^2+W(f_t)\right) \\
= & -|\nabla df_t|^2+\langle dd^*df_t, df_t\rangle - \text{Ric}(df_t, df_t)+\langle d\dot{f}_t, df_t\rangle \\
 &  +\langle \nabla W(f_t), d^*df_t\rangle+\langle \nabla W(f_t), \dot{f}_t\rangle  - \text{Hess}(W)_{AB}\langle df_t^A, df_t^B\rangle \\
= & -|\nabla df_t|^2-|\nabla W(f_t)|^2-\text{Ric}(df_t, df_t)-2\text{Hess}(W)_{AB}\langle df_t^A, df_t^B\rangle.
\end{split}
\end{equation}

Here $\dot{f}_t = \dsst \frac{\partial f_t}{\partial t}$ and we apply the equation $\dot{f}_t+d^*df_t+\nabla W(f_t)=0$ in the last step.

In Chen and Struwe (\cite{StruweSigma}), $W(f_t)=K\chi(\text{dist}(f_t, N))$. Here, $\chi(x)$ is a smooth function, $\chi(x)=x$ when $x<R_M$, the lower bound of the injectivity radius of $M$, and $\chi(x)$ is constant when $x\geq R_M$. $\text{dist}(f_t, N)$ is the distance function between the image of $f_t$ and $N$. Since $\chi(\text{dist}(f_t, N))$ is a bounded function, one can show that both $|\nabla W|^2$ and $\text{Hess}(W)_{AB}$ can be bounded by some constant depending on $W$, and therefore, we have the inequality
\begin{equation}
\left(\frac{\partial}{\partial t}+\Delta\right)e(f_t)\leq (C_1e(f_t)+C_2)e(f_t).
\end{equation}

Chen and Struwe \cite{StruweSigma} proved the following result
\begin{theorem}[Chen-Struwe \cite{StruweSigma}]
	Suppose $f:M\times \R_+\rightarrow N$ is a limit of a sequence $\{f_k\}$ with $\{f_k\}$ solving the parabolic harmonic equation \eqref{gradientflow} and with uniformly finite energy 
	$$ E(f_k)(t)\leq E_0<\infty,\quad  \quad \text{for $t>0$}, $$
	in the sense that $E(f)\leq E_0$ almost everywhere and $df_k\rightarrow df$ weakly in $L^2(Q)$ for any compact set $Q\in M\times \R_+$. Then $f$ solves the parabolic harmonic equation \eqref{gradientflow} in the classical sense and is regular on a dense open set whose complement $\Sigma$ has locally finite Hausdorff co-dimension $2$ with respect to the parabolic metric.  
\end{theorem}

This thesis began as an exploration of the sigma-model approach to harmonic maps in the context of Grassmannian target spaces. Grassmannians are of fundamental interest due to their role as classifying spaces. They are equipped with a wealth of explicitly computable algebraically defined potential energy functions, and they are positively curved. The latter feature means that we have no effective means for determining which maps
are homotopic to stable harmonic maps and which are not. Perhaps in such an explicit example, one might find new criteria guaranteeing the existence or nonexistence of stable harmonic maps.

It is convenient to identify the Grassmannian $\textrm{Gr}(k,\mathds{C}^l)$ of $k$ planes in $\mathds{C}^l$ with the space of hermitian involutions (rather than hermitian projections). With this identification, we note that the connected component of the space of rank $l$ hermitian matrices with nonzero determinant which contains $\textrm{Gr}(k,\mathds{C}^l)$ is homotopically equivalent to $\textrm{Gr}(k,\mathds{C}^l)$. Hence for topological purposes, it may be of interest to replace the target space $\textrm{Gr}(k,\mathds{C}^l)$ with this connected component in the
harmonic mapping problem.

In making this replacement, we are led to modify the usual sigma model approach where one chooses a sequence of potential energies which in some limit force a map to lie on the desired target space, with a single repulsive potential $W$ which simply prevents the map from entering a forbidden region - $\{f:\det(f) = 0\}$. In particular, we choose $W$ so that
$W(f) = \infty$ if and only if $\det(f)= 0$, and $W$ is smooth elsewhere.

In this new problem, issues of the singularity of a minimizing map $f$ become translated into possibly more accessible and explicit issues about singularities of $W(f)$. In fact, in order to show minimizers are smooth, one simply has to find conditions under which the potential energy remains bounded. We have not yet found such conditions. In these early explorations, we have been led first to sharpen the $\epsilon$-regularity
and monotonicity techniques which form the core of the regularity theory of harmonic maps and test them against explicit choices of algebraic potentials.

For example we consider the potential
$$W_b^L:=\textrm{trace}((f^{2L}+bI)^{-1})$$
where $I$ is the $l\times l$ identity matrix.  We would prefer to set $b=0$, but for many arguments, $b>0$ is convenient.
Thus we need effective estimates in order to prevent a minimizer from wandering out of the desired component. Of course, we do not expect such a result to hold without additional (as yet undiscovered) hypotheses. So instead, we first test our techniques by estimating the size of the bad
set.

We obtain the following theorems of bounds of $e_b^L(f)=\dsst \frac{1}{2}|df|^2+W_b^L(f)$:
\begin{theorem}[See Theorem 20]
	Suppose $f$ solves the elliptic equation 
	\begin{equation}\label{introWbL}
	 d^*df+W_b^L(f)=0 \qquad \textrm{with $E(f)$ finite}.
	 \end{equation} 
	
	There exists $\epsilon_0 = \dsst\frac{b^{\frac{1}{L}}}{2C}$, with $C$ depending only on $M, l, L$, such that if there exists $R<\rho$, and
	\begin{equation}
	\Phi_b^L(x_0, R) = \frac{1}{R^{m-2+p_0}}\int_{B_R(x_0)} (1+\nu(r, x))e_b^L(f)\, d\textup{vol} \leq \epsilon_0,
	\end{equation}
	then for any $\delta \leq \frac{3}{4}$, 
	\begin{equation}
	\sup_{B_{\delta R}(x_0)} e_b^L(f) \leq \frac{c(L, l, M)b^{\frac{1}{L}}}{(\delta R)^{2-p_0}}.
	\end{equation}
	Here, $\nu(r, x)$ and $p_0$ are defined as in Chapter $2$ \eqref{mu0nu0} and \eqref{p0}.
\end{theorem}

\begin{theorem}[See Theorem 21]
	Suppose $f_t$ solves the parabolic equation 
	\begin{equation}
	\frac{\partial f_t}{\partial t}+d^*df_t+W_b^L(f_t)=0 \qquad \textrm{with $f_0\in \textrm{Gr}(k, \mathds{C}^l)$ and $E(f_0)$ finite}.
	\end{equation}
	Here $f_t$ exists and is smooth for some $T<\infty$.
	
	There exists $\epsilon_0=\dsst \frac{b^{\frac{1}{L}}}{2C}$ with $C$ depending only on $M, l, L$, such that if there exists $R$ with $0<R<\sqrt{t_0-t}/2$ and
	\begin{equation}
	\Psi((x_0, t_0), R) = \int_{t_0-4R^2}^{t_0-R^2}\frac{1}{(t_0-t)^{\nu_0}}\int_M G_{(x_0, t_0)}(x, t)e_b^L(f_t)\varphi^2\, d\textup{vol}\, dt \leq \epsilon_0,
	\end{equation}
	then for any $\delta$ depending only on $M, R$, we have
	\begin{equation}
	\sup_{P_{\delta R}(z_0)} e_b^L(f) \leq \frac{c(L, l, M)b^{\frac{1}{L}}}{(\delta R)^{2-2\nu_0}}.
	\end{equation} 
	Here $z_0=(x_0,t_0)$, $\nu_0$ is the same ratio defined in \eqref{mu0nu0} and the theorem holds when $\nu_0=0$. 
\end{theorem}

When $L=1$, we also derive the absolute bound of $e_b(f)$ over the entire manifold $M$. Here we denote $\dsst \sup_M e_b(f)=e_M$.
\begin{equation}
e_M \leq C(M)E_b(f)^{\frac{2m}{2m-2}}b^{-\frac{2m-1}{m-1}}.
\end{equation}
(See \eqref{eMbound})

The upper bound of the measure of the bad set defined as
\begin{equation}\label{introSigmabL}
 \Sigma_b^L = \left\{x\in M|e_b^L(f)\geq \frac{1}{b} \right\}
 \end{equation}
 is given by the following theorem
\begin{theorem}[See Theorem 22]
	Fix a constant $b>0$. Let $f$ be the solution to \eqref{introWbL}. Define the set $\Sigma_b^L$ as in \eqref{introSigmabL}. The Hausdorff $\left(\frac{2}{L+1}+(m-2)\right)$-measure of $\Sigma_b^L$ at scale $b^{\frac{1}{2}+\frac{1}{2L}}$ (defined in \eqref{hausdorff}) is bounded and independent of $b$.
\end{theorem}

%% file: Chapter2/Chapter2.tex
\def\dsst{\displaystyle}
\def\R{\mathds{R}}
\chapter{Monotonicity Formulas}
\label{chap:price}

When an energy density $e$ satisfies a partial differential inequality of the form \\$\Delta e \leq Ce$, Moser iteration allows one to obtain pointwise estimates for $e$ in terms of the integral of $e$ over balls. When studying nonlinear equations, the preceding equality is often replaced by an inequality of the form  $\Delta e \leq Ce^p,$ for some $p>1$. In such situations, Moser iteration must be supplemented by $\epsilon$-regularity arguments. A crucial component of such arguments are monotonicity relations for rescaled local energies. In this section, we derive such monotonicity relations in the context of reaction diffusion problems.

Our derivation is more flexible than the standard method (e.g., \cite{price1983monotonicity}); moreover, it clarifies how to find optimal monotonicity relations.

We first introduce some basic notation. Denote by $|x-y|$ the geodesic distance between $x$ and $y$ in $M$.
For $x_0\in M$ and $R>0$ let
\begin{equation}
	B_R(x_0) = \left\{x\in M| |x-x_0|\leq R \right\},
\end{equation}
and
\begin{equation}
S_R(x_0) = \left\{x\in M||x-x_0|=R \right\}.
\end{equation}

For the parabolic case we need additional notation. Let $z=(x,t) \in M\times \R$. For a distinguished point $z_0=(x_0,t_0)$ and $R>0$, let
\begin{equation}
P_R(z_0) = \left\{z=(x,t)| |x-x_0|\leq R,\ |t-t_0|\leq R^2 \right\},
\end{equation}
\begin{equation}
S_R(z_0) = \left\{z=(x,t) | |x-x_0|=R,\ t=t_0 \right\},
\end{equation}
and
\begin{equation}
T_R(t_0) = \left\{z=(x,t) | t_0-4R^2\leq t\leq t-R^2 \right\}.
\end{equation}

Let $M$ be a compact Riemannian manifold. Consider a function $f:M\rightarrow \R^k.$ Let $W:\R^k\rightarrow \R$ be a function, which we will initially take to be smooth, but will later allow to be singular.  Define the energy density
\begin{equation}
e(f):= \frac{1}{2}|df|^2+ W(f).
\end{equation}
We will be interested in whether the kinetic term $\frac{1}{2}|df|^2$ or the potential term $W(f)$ dominates the energy density. Therefore we define the following ratios over geodesic spheres $S_R(x)\subset M$:

\begin{equation}\label{mu}
\mu(R, x) = \frac{\int_{S_R(x)}|f_r|^2\, d\sigma}{\int_{S_R(x)}\frac{1}{2}|df|^2+W(f)\, d\sigma},
\end{equation}
and
\begin{equation}\label{nu}
\nu(R,x) = \frac{\int_{S_R(x)}W(f)\, d\sigma}{\int_{S_R(x)}\frac{1}{2}|df|^2+W(f)\, d\sigma}.
\end{equation}

We will see that nontrivial lower bounds for $\mu(R,x)$ and $\nu(R,x)$ allow us to sharpen estimates. We will need such lower bounds on balls of small, but not arbitrarily small radius. Hence we will fix a scale $\rho >0$ and define
\begin{equation}
\mu(x) =\inf_{R\leq \rho}\mu(R,x),\qquad \nu(x) =\inf_{R\leq \rho}\nu(R,x).
\end{equation}
When necessary for clarity, we will write instead $\mu_\rho(x)$ to make the scale explicit.

We also set
\begin{equation}\label{mu0nu0}
\mu_{0,S}:= \inf_{x\in S}\mu(x), \qquad \nu_{0, S}:=\inf_{x\in S}\nu(x).
\end{equation}
When $S$ is clear from context, we will simply write $\mu_0$ and $\nu_0$.

It is obvious that $0\leq \nu(x), \nu_0\leq 1$. From \cite{di2017price},
\begin{equation}
\lim_{R\rightarrow 0}\mu(R,x) = \frac{2}{m}.
\end{equation}
Therefore, we have $\mu(x), \mu_0\leq \dsst \frac{2}{m}$. 

With these preliminaries, we now turn to the derivation of monotonicity relations.

\section{Monotonicity Formula for Elliptic Equations}

Let $R_M$ denote the lower bound of the injectivity radius of $M$.

\begin{lemma}\label{monoelliptic}
	Suppose $f(x):M\rightarrow \R^k$ is a smooth solution to the differential equation
	\begin{equation}
	d^*df+W(f) = 0.
	\end{equation}
  Fix a scale $\rho<R_M$. Let $\dsst p_0=\frac{2\nu_{0, B_R(x_0)}+(m-2)\mu_{0, B_R(x_0)}}{1-\mu_{0, B_R(x_0)}}$. Define 
	\begin{equation}
	\Phi(R, x) = e^{CR}R^{2-m-p_0}\int_{B_R(x)}\left(\frac{1}{2}|df|^2+\frac{m}{m-2}W(f)\right)\, d\textup{vol},
	\end{equation}
	where $m=\dim(M)$, $C$ is a constant which only depends on the geometry of $M$. Then $\Phi(R,x)$ is non-decreasing for any $x\in M$, $R<\rho$.
\end{lemma}

\begin{proof}
	Without loss of generality, we assume $x=0$, the origin of a geodesic coordinate system $\{x^j\}_{j=1}^m$, and denote $B_R(0)=B_R$. For a 1-form $\omega$, let $e(\omega)$ denote exterior multiplication on the left by $\omega$, and let
	$e^*(\omega)$ denote its adjoint. Let $r$ denote the radial function of the geodesic coordinates. Let us consider
	\begin{equation} \label{ellipquan}
	\int_{B_R} \langle\{d, e^*(rdr)\}df, df\rangle \, \textrm{vol}.
	\end{equation}
	
	We compute this quantity in two ways. First, by integration by parts, we have
	\begin{equation}\label{quantity1}
	\begin{split}
	& \int_{B_R}\langle \{d, e^*(rdr)\}df, df\rangle \, d\textrm{vol}\\
	 = &\int_{B_R}\langle d(e^*(rdr)df), df\rangle\, d\textrm{vol}\\
	= & \int_{B_R} d(e^*(rdr)df\wedge \ast df) + \int_{B_R}\langle e^*(rdr)df, d^*df\rangle \, d\textrm{vol}\\
	= & \int_{S_R} r|e^*(dr)df|^2\, d\sigma +\int_{B_R} \langle e^*(rdr)df, d^*df\rangle\, d\textrm{vol}
	\end{split} 
	\end{equation}
	
	Define
	\begin{equation}
	Q= e(dx^j)e^*(\nabla_j dr).
	\end{equation}
	
	$Q$ is a natural extension of the second fundamental form of the geodesic sphere to an endomorphism of differential forms. We can compute the action of $Q$ on $k$-forms as follows:
	\begin{equation}\label{Qformula}
	Q =e(dx^j)e^*(\nabla_j \frac{x_j}{r}dr) =  \frac{e(dx^i)e^*(dx^i)}{r}-\frac{e(dr)e^*(dr)}{r} = \frac{k-e(dr)e^*(dr)}{r}+A(x),
	\end{equation}
	where $A(x)=O(r)$ arises from connection terms and $g_{ij} - \delta_{ij}$.
	
	Next, we use \eqref{Qformula} to compute \eqref{ellipquan},
	\begin{equation}\label{quantity2}
	\begin{split}
	& \langle \{d, e^*(rdr)\}df, df\rangle_{L^2(B_R)} \\
	= & \langle \left\{e(dx^j)\nabla_j, e^*(rdr) \right\} df, df\rangle_{L^2(B_R)}\\
	= & \langle \left(e(dx^j)[\nabla_j, e^*(rdr)]+\left\{e(dx^j), e^*(rdr) \right\}\nabla_j\right)df, df\rangle_{L^2(B_R)}\\
	= & \langle e(dx^j)e^*(\nabla_j(rdr))df, df\rangle_{L^2(B_R)}+\int_{B_R} r\frac{\partial}{\partial r}\left(\frac{1}{2}|df|^2\right)\, d\textrm{vol}\\
	= &  \int_{B_R} r\frac{\partial}{\partial r}\left(\frac{1}{2}|df|^2\right)\, d\textrm{vol} +\int_{B_R} \langle e(dr)e^*(dr)df, df\rangle\, d\textrm{vol} +\int_{B_R} r\langle Qdf,df\rangle\, d\textrm{vol}\\
	= & \int_{B_R}r\left(\{d, e^*(dr)\}-Q\right)\left(\frac{1}{2}|df|^2\, dvol\right)+\int_{B_R}|df|^2\, d\textrm{vol} + \int_{B_R} r\langle A(x)df, df\rangle\, d\textrm{vol}\\ 
	= & \int_{B_R} r\left(d\iota_r+\iota_r d\right)\left(\frac{1}{2}|df|^2\, d\textrm{vol}\right)-\int_{B_R}rQ\left(\frac{1}{2}|df|^2\, d\textrm{vol}\right)+\int_{B_R}|df|^2\, d\textrm{vol}\\
	& + \int_{B_R} r\langle A(x)df, df\rangle\, d\textrm{vol}\\ 
	= & \int_{B_R} d\left(r\iota_r\frac{1}{2}|df|^2\, d\textrm{vol}\right)-\int_{B_R}\frac{1}{2}|df|^2\, d\textrm{vol}-\int_{B_R}r\frac{1}{2}|df|^2Q\left(d\textrm{vol}\right)\\
	& + \int_{B_R}|df|^2\, d\textrm{vol} + \int_{B_R} r\langle A(x)df, df\rangle\, d\textrm{vol}\\ 
	= & \int_{S_R} R\left(\frac{1}{2}|df|^2\right)\, d\sigma-\int_{B_R}\frac{1}{2}|df|^2\, d\textrm{vol}-\int_{B_R}r\frac{1}{2}|df|^2Q\left(d\textrm{vol}\right)\\
	& \int_{B_R}|df|^2\, d\textrm{vol} + \int_{B_R} r\langle A(x)df, df\rangle\, d\textrm{vol}.
	\end{split}
	\end{equation}
	
	We first compute the anticommutator in \eqref{quantity1} in terms of $Q$ and then integrate.
	\begin{equation}\label{QonVol}
	Q(d\textrm{vol}) = H d\textrm{vol} = \left(\frac{m-1}{r}+A(x)\right)\,d\textrm{vol},
	\end{equation}
	where, under an orthonormal frame, 
	\begin{equation}
	H = \sum_k^m \langle \nabla_k\left(\partial r\right), \partial_k\rangle.
	\end{equation}
	
	Applying \eqref{QonVol} on \eqref{quantity2}, we get
	\begin{equation}
	\begin{split}
	& \langle {d, e^*(rdr)}df, df\rangle_{L^2(B_R)}\\
	= & \int_{S_R}R\frac{1}{2}|df|^2\, d\sigma - \int_{B_R}\left(\frac{m}{2}+\frac{rA(x)}{2}\right)|df|^2\, d\textrm{vol} + \int_{B_R}(|df|^2+r\langle A(x)df, df\rangle)\, d\textrm{vol}\\
	= & \int_{B_R}\left(1- \frac{m}{2}+B(x)\right)|df|^2 \, d\textrm{vol} +\int_{S_R}R\frac{1}{2}|df|^2\, d\sigma.
	\end{split}
	\end{equation}
	where $B(x)|df|^2=\dsst -\frac{rA(x)}{2}|df|^2+r\langle A(x)df, df\rangle = O(r^2)|df|^2$.
	
	Combining \eqref{quantity1} and \eqref{quantity2}, we have
	\begin{equation}\label{quancombined}
	\begin{split}
	& \int_{B_R}\left(1-\frac{m}{2}+B(x)\right)|df|^2 \, d\textrm{vol} \\
	= & \int_{B_R} \langle e^*(rdr)df, d^*df\rangle \, d\textrm{vol} + \int_{S_R} R\left[|e^*(dr)df|^2-\frac{1}{2}|df|^2\right]\, d\sigma. 
	\end{split}
	\end{equation}
	
	Since $f$ satisfies the elliptic equation
	$$ d^*df+\nabla W(f) = 0, $$
	we can replace $d^*df$ by $-\nabla W(f)$ on the right-hand side of \eqref{quancombined} to get
	\begin{equation}
	\begin{split}
	& \int_{B_R} \langle e^*(rdr)df, d^*df\rangle\, d\textrm{vol} \\
	= & -\int_{B_R} \langle e^*(rdr)df, \nabla W(f)\rangle\, d\textrm{vol} \\
	= & - \int_{B_R} r\frac{\partial}{\partial r}W\, d\textrm{vol} =  -\int_{S_R} RW\, d\sigma + \int_{B_R} \left(m+\hat{B}(x)\right) W\, d\textrm{vol}.
	\end{split}
	\end{equation}
	where $\hat{B}(x)= O(r^2)$ consists of geometric terms. In the last step we use Stokes' Theorem.
	
	Therefore
	\begin{equation}\label{ellipticequality}
	\begin{split}
	& \int_{B_R}\left[\left(1-\frac{m}{2}+B(x)\right)|df|^2 -\left(m+\hat{B}(x)\right)W\right]\, d\textrm{vol} \\
	= & \int_{S_R}R\left[|e^*(dr)df|^2-\frac{1}{2}|df|^2-W\right]\, d\sigma.
	\end{split}
	\end{equation}
	
	
	
	
	Using the ratios $\mu$ and $\nu$ defined in \eqref{mu} and \eqref{nu}, \eqref{ellipticequality} can be rewritten
	\begin{equation}\label{ellipticbeforephi}
	\begin{split}
	& \int_{B_R}\left(m-2+2\nu(1+\hat{B}(x))-2B(x)\right) \left(\frac{1}{2}|df|^2+W\right)d\textrm{vol}\\
	= &   \int_{S_R}R\left(1-\mu\right)\left(\frac{1}{2}|df|^2+W\right)\, d\sigma.
	\end{split}
	\end{equation}
	
	Let $\phi(R)$ be a function to be determined later. Multiply \eqref{ellipticbeforephi} by $\phi'(R)$ and integrate from $\sigma$ to $\tau$ to obtain
	\begin{equation}\label{LHS}
	\begin{split} 
	& \phi(\tau)\int_{B_\tau}(m-2+2\nu)\left(\frac{1}{2}|df|^2+W\right)\, d\textrm{vol} -\phi(\sigma)\int_{B_\sigma}(m-2+2\nu)\left(\frac{1}{2}|df|^2+W\right)\, d\textrm{vol}\\
	&-\int_{B_\tau\setminus B_\sigma}\phi(r)(m-2+2\nu)\left(\frac{1}{2}|df|^2+W\right)\, d\textrm{vol}\\
	& + \int_{\sigma}^{\tau}\phi'(R)\int_{B_R}D(x)\left(\frac{1}{2}|df|^2+W\right)\, d\textrm{vol}\, dR\\
	= & \int_{B_\tau\setminus B_\sigma}\phi'(r)r\left(1-\mu\right)\left(\frac{1}{2}|df|^2+W\right)\, d\textrm{vol}, 
	\end{split}
	\end{equation}
	where $D(x) = 2\nu \hat{B}(x)-2B(x)=O(r^2)$.
	
	Rearranging the order, we get
	\begin{equation}
	\begin{split} 
	& \phi(\tau)\int_{B_\tau}(m-2+2\nu)\left(\frac{1}{2}|df|^2+W\right)\, d\textrm{vol} -\phi(\sigma)\int_{B_\sigma}(m-2+2\nu)\left(\frac{1}{2}|df|^2+W\right)\, d\textrm{vol}\\
	= &\int_{B_\tau\setminus B_\sigma}\left[\phi(r)\left(m-2+2\nu\right)+\phi'(r)r(1-\mu)\right]\left(\frac{1}{2}|df|^2+W\right)\, d\textrm{vol}\\
	& + \int_{\sigma}^{\tau}\phi'(R)\int_{B_R}D(x)\left(\frac{1}{2}|df|^2+W\right)\, d\textrm{vol}\, dR.
	\end{split}
	\end{equation}
	
	Now pick $\phi(r)$ to be
	\begin{equation}
	\phi(\tau) = \phi(\sigma)\exp\left(\int_{\sigma}^{\tau}\frac{2-m-2\nu}{r(1-\mu)}\, dr\right)=\phi(\sigma)\exp\left(\int_{\sigma}^{\tau}\frac{2-m}{r}-\frac{2\nu+(m-2)\mu}{r(1-\mu)}\, dr\right).
	\end{equation}
	Then
	\begin{equation}
	\phi(r)\left(m-2+2\nu\right)+\phi'(r)r(1-\mu) = 0 .
	\end{equation}
	Therefore,
	\begin{equation}\label{finalelliptic}
	\begin{split} 
	& \phi(\tau)\int_{B_\tau}(m-2+2\nu)\left(\frac{1}{2}|df|^2+W\right)\, d\textrm{vol} -\phi(\sigma)\int_{B_\sigma}(m-2+2\nu)\left(\frac{1}{2}|df|^2+W\right)\, d\textrm{vol}\\
	= & \int_{\sigma}^{\tau}\phi'(R)\int_{B_R}D(x)\left(\frac{1}{2}|df|^2+W\right)\, d\textrm{vol}\, dR.
	\end{split}
	\end{equation}
	
	Let
	\begin{equation}\label{p0}
	p_0=\frac{2\nu_0+(m-2)\mu_0}{1-\mu_0} .
	\end{equation}
	When $m>2$, 
	\begin{equation}
	p_0=\frac{2\nu_0+(m-2)\mu_0}{1-\mu_0} \geq 0.
	\end{equation}
	We have
	\begin{equation}\label{phicompare}
	\phi(\tau)\leq \phi(\sigma)\left(\frac{\tau}{\sigma}\right)^{2-m-p_0}.
	\end{equation}
	
	Since $D(x)= O(r^2)$ and is bounded, on $B_R$ we can assume $|D(x)|\leq CR$, where $C\geq0$. 
	
	Hence, 
	\begin{equation}
	\begin{split}
	\left|\int_{\sigma}^{\tau}\phi'(R)\int_{B_R}D(x)e\,d\textrm{vol}\, dR\right|\leq &\int_{\sigma}^{\tau} \left|\phi'(R)\right|\int_{B_R}|D(x)|e\, d\textrm{vol} \, dR \\
	\leq & C\int_{\sigma}^{\tau} \left|\phi'(R)\right|R\int_{B_R}e\, d\textrm{vol}\, dR. 
	\end{split}
	\end{equation}
	
	Since $\phi'(R)\leq 0$, using
	\begin{equation}
	|\phi'(R)|R=-\phi'(R)R=\frac{(m-2+2\nu)}{1-\mu}\phi(R),
	\end{equation} 
	we get
	\begin{equation}\label{estimateforerrorterm}
	\begin{split}
	&\int_{\sigma}^{\tau}\phi'(R)\int_{B_R}D(x)e\, d\textrm{vol}\, dR\\
	 \geq & -C\int_{\sigma}^{\tau}\phi(R)\int_{B_R}\frac{m-2+2\nu}{1-\mu}e\, d\textrm{vol}\, dR \\
	\geq & -2C\int_{\sigma}^{\tau}\phi(R)\int_{B_R}(m-2+2\nu)e\, d\textrm{vol}\, dR\\
	= & -2C\int_{\sigma}^{\tau}\phi(R)\int_{B_R}(m-2+2\nu)\left(\frac{1}{2}|df|^2+W(f)\right)\, d\textrm{vol}\, dR.
	\end{split}
	\end{equation}
	
	Define 
	$$ h(R)=\phi(R)\int_{B_R}(m-2+2\nu)\left(\frac{1}{2}|df|^2+W(f)\right)\, d\textrm{vol}. $$
	
	We combine \eqref{finalelliptic} and \eqref{estimateforerrorterm} to get
	\begin{equation}
	h(\tau)-h(\sigma)\geq -2C\int_{\sigma}^{\tau} h(R)\, dR. 
	\end{equation}
	
	Dividing both sides by $\tau-\sigma$ and taking the limit $\tau\rightarrow \sigma$, this becomes 
	$$ h'(\sigma)\geq -2Ch(\sigma). $$
	
	We solve $h$ to get
	$$ e^{2C\tau}h(\tau)\geq e^{2C\sigma}h(\sigma), $$
	which is
	\begin{equation}\label{errorterm}
	e^{2C\tau}\phi(\tau)\int_{B_\tau}(m-2+2\nu)e\, d\textrm{vol} \geq e^{2C\sigma}\phi(\sigma)\int_{B_\sigma}(m-2+2\nu)e\, d\textrm{vol}.
	\end{equation}	


	

	Combining \eqref{finalelliptic}, \eqref{phicompare}, and \eqref{errorterm}, we get
	
	\begin{equation}
	\begin{split}
	&e^{2C\tau}\tau^{2-m-p_0}\int_{B_\tau}(m-2+2\nu)\left(\frac{1}{2}|df|^2+W\right)\, d\textrm{vol} \\
\geq & e^{2C\sigma}\sigma^{2-m-p_0}\int_{B_\sigma}(m-2+2\nu)\left(\frac{1}{2}|df|^2+W\right)\, d\textrm{vol}.
	\end{split}
	\end{equation}
	
	That is,
	\begin{equation}
	\begin{split}
	& e^{2C\tau}\tau^{2-m-p_0}\int_{B_\tau}\left(\frac{1}{2}|df|^2+\frac{m}{m-2}W\right)\, d\textrm{vol} \\
	\geq&  e^{2C\sigma}\sigma^{2-m-p_0}\int_{B_\sigma}\left(\frac{1}{2}|df|^2+\frac{m}{m-2}W\right)\, d\textrm{vol}.
	\end{split}
	\end{equation}
	
	When $M$ is a Euclidean space, $C=0$ since $D(r)=0$.
\end{proof}

\section{Monotonicity Formula for Parabolic Equations}

\begin{lemma}
	Suppose $f(t,x)=f_t(x):M\times I\rightarrow \R^k$ is a regular solution to the differential equation
	\begin{equation}\label{paraequa}
	\frac{\partial f_t}{\partial t}+d^*df_t+W(f_t) = 0,
	\end{equation}
	where $I$ is a subset of the interval when $f_t(x)$ exists in short time. Let $t_0\in I$, and $0<R<R_0<R_M$ such that $[t_0-4R^2,t_0-R^2]\subset I$.\\
	
	Define
	\begin{equation}
	\phi(R, x_0, t_0) = \frac{1}{R^{2-2\nu_0}}\int_M G_{(x_0,t_0)}(x,t)e(f)\varphi^2(x)\, d\textup{vol},
	\end{equation}
	and
	\begin{equation}\label{Psi}
	\begin{split}
	\Psi(R,x_0,t_0) = & \int_{t_0-4R^2}^{t_0-R^2}\frac{\phi(\sqrt{t_0-t})}{t_0-t}\, dt\\
	= &\int_{t_0-4R^2}^{t_0-R^2}\frac{1}{(t_0-t)^{\nu_0}}\int_M G_{(x_0,t_0)}(x,t)e(f)\varphi^2(x)\, d\textup{vol}\, dt,
	\end{split}
	\end{equation}
	where $\nu_0$ is the same ratio defined in \eqref{mu0nu0}, 
	\begin{equation}
	G_{(x_0, t_0)}(x,t)=\frac{1}{(4\pi(t_0-t))^{\frac{m}{2}}}\exp\left(-\frac{|x-x_0|^2}{4\pi(t_0-t)}\right)
	\end{equation}
	is the Euclidean backward Gaussian, and $\varphi(x)=\varphi(|x-x_0|)$ is a compactly supported radial function on $B_{R_M}(x_0)$  such that $\varphi(R)=1$ for $R<R_M/2$ and $|\varphi'(R)|\leq \frac{4}{R_M}$.  
	
	Then we have the following monotonicity formulas for $\phi(R)$ and $\Psi(R)$:
	\begin{equation}\label{phiinequa}
	\phi(R)\leq \exp\left(c(R_0^2-R^2)\right)\phi(R_0)+\hat{C}E(f_0(x))(R_0^2-R^2),
	\end{equation}
	and
	\begin{equation}\label{psiinequa}
	\Psi(R)\leq \exp(c(R_0^2-R^2))\Psi(R_0)+\hat{C}E(f_0(x))(R_0^2-R^2),
	\end{equation}
	for any $0<R<R_0<\sqrt{t_0-t}/2$. Here constants $c$, $\hat{C}$ only depend on the geometry of $M$, and $m=\dim(M)$. $E$ is the energy of $f_t$ over $M$:
	\begin{equation}
	E(f_t) = \int_M \left(\frac{1}{2}|df_t|^2+W(f_t)\right)\, d\textup{vol}. 
	\end{equation}
\end{lemma}

It is clear that when $\nu(R,x)=0$, which implies $\nu_0=0$, the above results \eqref{phiinequa} and \eqref{psiinequa} reduce to results of Struwe \cite[Lemma 4.2]{StruweSigma}.
\begin{proof}
	Without loss of generality, we assume $x_0=0$. Let $\varphi(r)$ be the radial cut-off function in $B_R(0)$ for some $R<R_M$. Consider
	\begin{equation}\label{parapositive}
	\int_M \langle \{d, e^*(dY)\}df, \varphi^2 df\rangle\, d\textrm{vol},
	\end{equation}
	where $Y$ is another radial function.
	
	Our goal is to find a $Y$ so that \eqref{parapositive} produces new monotonicity result. We calculate this quantity in two ways. On one hand, 
	\begin{equation}
	\begin{split}
	& \int_M \langle \{d, e^*(dY)\}df, \varphi^2df\rangle\, d\textrm{vol} \\
	= & \int_M\langle \left\{e(dx^j)\nabla_j, e^*(dY) \right\} df, \varphi^2df\rangle\, d\textrm{vol}\\
	= & \int_M \langle \left(e(dx^j)[\nabla_j, e^*(dY)]+\left\{e(dx^j), e^*(dY) \right\}\nabla_j\right)df, \varphi^2df\rangle\, d\textrm{vol}\\
	= & \int_M \langle e(dx^j)\nabla_j\left(Y_re^*(dr)\right)df, \varphi^2 df\rangle\, d\textrm{vol}+\int_M Y_r\langle \nabla_r(df), \varphi^2df\rangle\, d\textrm{vol} \\
	= & \int_M \nabla_j(Y_r)\varphi^2\langle  e(dx^j)e^*(dr)df, df\rangle\, d\textrm{vol} + \int_M Y_r\varphi^2\langle e(dx^j)e^*(\nabla_j dr)df, df\rangle\, d\textrm{vol}\\
	&  + \int_M Y_r \frac{\partial}{\partial r}\left(\frac{1}{2}|df|^2\right)\varphi^2\, d\textrm{vol} \\
	= & \int_M Y_{rr}|e^*(dr)df|^2 \varphi^2\, d\textrm{vol} +\int_M Y_r\langle Qdf, df\rangle \varphi^2\, d\textrm{vol} + \int_M Y_r \frac{\partial}{\partial r}\left(\frac{1}{2}|df|^2\right)\varphi^2\, d\textrm{vol}\\
	= & \text{I} + \text{II} + \text{III}.
	\end{split}
	\end{equation} 
	
	For the second term $\text{II}$, we use
	\begin{equation}
	Q=e(dx^j)e^*(\nabla_j dr)=\frac{1-e(dr)e^*(dr)}{r}+A(x),
	\end{equation}
	where $A(x)= O(r)$, to replace $Q$:
	\begin{equation}
	\text{II} =  \int_M Y_r\left(\frac{1}{r}(|df^2|-|e^*(dr)df|^2) + \langle A(x)df,df\rangle\right)\varphi^2\, d\textrm{vol}.
	\end{equation}
	
	Integration by parts in the third term $\text{III}$, we have
	\begin{equation}
	\begin{split}
	\text{III} = & \int_M Y_r\varphi^2\frac{\partial}{\partial r}\left(\frac{1}{2}|df|^2\right)\, d\textrm{vol} = \int_M \langle dY, d\left(\varphi^2\frac{1}{2}|df|^2\right)\rangle \, dx - \int_M Y_r|df|^2\varphi\varphi'\, d\textrm{vol}\\
	= & \int_M d^*dY\left(\frac{1}{2}|df|^2\right)\varphi^2\, d\textrm{vol} - \int_M Y_r|df|^2\varphi\varphi'\, d\textrm{vol}.
	\end{split}
	\end{equation} 
	
	Therefore,
	\begin{equation}\label{paraquantity1}
	\begin{split}
	& \int_M \langle \{d, e^*(dY)\}df, \varphi^2 df\rangle\, d\textrm{vol} \\
	= & \text{I} + \text{II} + \text{III} \\
	= & \int_M Y_{rr}|e^*(dr)df|^2\varphi^2\, d\textrm{vol} + \int_M Y_r\left(\frac{1}{r}(|df|^2-|e^*(dr)df|^2)+\langle A(x)df, df\rangle\right)\varphi^2\, d\textrm{vol}\\
	& +\int_M d^*dY\left(\frac{1}{2}|df|^2\right)\varphi^2\, d\textrm{vol} - \int_M Y_r|df|^2\varphi\varphi'\, d\textrm{vol}.
	\end{split}
	\end{equation}

	On the other hand, we may apply integration by parts at the very beginning and get another expression about $\langle \{d, e^*(dY)\}df, \varphi^2df\rangle_{L^2(M)}$,
	\begin{equation}
	\begin{split}
	& \int_M \langle \{d, e^*(dY)\}df, \varphi^2df\rangle\, d\textrm{vol} \\
	= & \int_M \langle e^*(dY)df, d^*\left(\varphi^2df\right) \rangle\, d\textrm{vol}\\
	= & \int_M \langle e^*(dY)df, d^*df\rangle \varphi^2\, d\textrm{vol} - \int_M \langle e^*(dY)df, 2\varphi \varphi'(r)e^*(dr)df\rangle\, d\textrm{vol}.
	\end{split}
	\end{equation}
	
	For the first term, we replace $d^*df$ using \eqref{paraequa} to get
	\begin{equation}
	\begin{split}
	& \int_M \langle e^*(dY)df, d^*df\rangle\varphi^2\, d\textrm{vol} \\
	= & -\int_M Y_r\langle e^*(dr)df, f_t+\nabla W\rangle\varphi^2\, d\textrm{vol}\\
	= & -\int_M \langle df, e(dY)f_t\rangle \varphi^2\, d\textrm{vol} - \int_M Y_r e^*(dr)dW \varphi^2\, d\textrm{vol}\\
	= & -\int_M \langle df, d(Yf_t)\rangle\varphi^2\, d\textrm{vol} + \int_M \langle df, Yd(f_t)\rangle \varphi^2\, d\textrm{vol} -\int_M \langle dY, d(\varphi^2 W)\rangle\, d\textrm{vol}\\
	&  + \int_M Y_rW2\varphi\varphi'\, d\textrm{vol}\\
	= & -\int_M \langle df, d(Yf_t)\rangle\varphi^2\, d\textrm{vol} -\int_M \langle df, Yd(f_t)\rangle\varphi^2\, d\textrm{vol} +2\int_M \langle df, Y(df_t)\rangle\varphi^2\, d\textrm{vol}\\
	&- \int_M d^*dY W\varphi^2\, dx+\int_M Y_rW2\varphi\varphi'\, d\textrm{vol}\\
	= & -\int_M Y\frac{d}{dt}\left(\frac{1}{2}|df|^2\right)\varphi^2\, d\textrm{vol} -2\int_M \langle df, [d, Y]f_t\rangle\varphi^2\, d\textrm{vol} +\int_M \langle df, d(Yf_t\varphi^2)\rangle\, d\textrm{vol}\\
	& - \int_M Y\langle f_r, f_t\rangle 2\varphi\varphi'\, d\textrm{vol} - \int_M d^*dYW\varphi^2\, d\textrm{vol} + \int_M Y_rW2\varphi\varphi'\, d\textrm{vol}\\
	= & -\int_M Y\frac{d}{dt}\left(\frac{1}{2}|df|^2+W\right)\varphi^2\, d\textrm{vol} - 2\int_M Y_r\langle f_r, f_t\rangle \varphi^2\, d\textrm{vol} + \int_M Y\langle d^*df, f_t\rangle\varphi^2\, d\textrm{vol}\\ 
	 & -\int_M Y\langle f_r, f_t\rangle 2\varphi\varphi'\, d\textrm{vol}  + \int_M Y\langle \nabla W, f_t\rangle\varphi^2\, d\textrm{vol}   - \int_M d^*dYW\varphi^2\, dx \\
	 & + \int_M Y_rW2\varphi\varphi'\, d\textrm{vol}\\
	= & -\frac{d}{dt}\int_M Y\left(\frac{1}{2}|df|^2+W\right)\varphi^2\, d\textrm{vol} + \int_M Y_t\left(\frac{1}{2}|df|^2+W\right)\varphi^2\, d\textrm{vol} \\
	& -\int_Md^*dYW\varphi^2\, d\textrm{vol}-\int_M Y|f_t|^2\varphi^2\, d\textrm{vol}  -2\int_M Y_r\langle f_r, f_t\rangle\varphi^2\, d\textrm{vol}\\
	& - \int_M Y\langle f_r, f_t\rangle 2\varphi\varphi'\, d\textrm{vol} + \int_M Y_rW2\varphi\varphi'\, d\textrm{vol}.
	\end{split}
	\end{equation}
	In the last step, we use the equation \eqref{paraequa} again.
	
	Therefore,
	\begin{equation}\label{paraquantity2}
	\begin{split}
	& \int_M \langle \{d, e^*(dY)\}df, \varphi^2 df\rangle\, d\textrm{vol} \\
	= &  -\frac{d}{dt}\int_M Y\left(\frac{1}{2}|df|^2+W\right)\varphi^2\, d\textrm{vol} + \int_M Y_t\left(\frac{1}{2}|df|^2+W\right)\varphi^2\, d\textrm{vol} \\
	& -\int_Md^*dYW\varphi^2\, d\textrm{vol}-\int_M Y|f_t|^2\varphi^2\, d\textrm{vol}  -2\int_M Y_r\langle f_r, f_t\rangle\varphi^2\, d\textrm{vol} \\
	& - \int_M Y\langle f_r, f_t\rangle 2\varphi\varphi'\, d\textrm{vol}+ \int_M Y_rW2\varphi\varphi'\, d\textrm{vol} -\int_M Y_r|f_r|^22\varphi\varphi'\, d\textrm{vol}.
	\end{split}
	\end{equation}

	Equating \eqref{paraquantity1} and \eqref{paraquantity2}, we get
	\begin{equation}
	\begin{split}
	0 = & \frac{d}{dt}\int_M Y\left(\frac{1}{2}|df|^2+W\right)\varphi^2\, d\textrm{vol} + \int_M \left(d^*dY-Y_t\right)\left(\frac{1}{2}|df|^2+W\right)\varphi^2\, d\textrm{vol}\\
	&+ \int_M Y|f_t|^2\varphi^2\, d\textrm{vol}+2\int_M Y_r\langle f_r, f_t\rangle\varphi^2\, d\textrm{vol}  + \int_M\left(Y_{rr}-\frac{Y_r}{r}\right)|f_r|^2\varphi^2\, d\textrm{vol}\\
	& +\int_M \left(\frac{Y_r}{r}|df|^2+Y_r\langle A(x)df,df\rangle\right)\varphi^2\, d\textrm{vol}
	 +\int_M Y\langle f_r, f_t\rangle 2\varphi\varphi'\, d\textrm{vol}\\
	 & +\int_M Y_r|f_r|^22\varphi\varphi'\, d\textrm{vol}-\int_M Y_r\left(\frac{1}{2}|df|^2+W\right)2\varphi\varphi'\, d\textrm{vol}\\
	=  & \frac{d}{dt}\int_M Y\left(\frac{1}{2}|df|^2+W\right)\varphi^2\, d\textrm{vol}+ \int_M \left(d^*dY-Y_t+\frac{2Y_r}{r}\right)\left(\frac{1}{2}|df|^2+W\right)\varphi^2\, d\textrm{vol}\\
	& + \int_M Y\left|f_t+\frac{Y_r}{Y}f_r\right|^2\varphi^2\, d\textrm{vol}  + \int_M \left(Y_{rr}-\frac{Y_r}{r}-\frac{Y_r^2}{Y}\right)|f_r|^2\varphi^2\, d\textrm{vol} \\
	& -\int_M \frac{Y_r}{r}2W\varphi^2\, d\textrm{vol} +\int_M Y\langle f_r, f_t+\frac{Y_r}{Y}f_r\rangle 2\varphi\varphi'\, d\textrm{vol}\\
	& -\int_M Y_r\left(\frac{1}{2}|df|^2+W\right)2\varphi\varphi'\, d\textrm{vol}+\int_M Y_r\langle A(x)df, df\rangle\varphi^2\, d\textrm{vol}.
	\end{split}
	\end{equation}
	
	Recall $\nu(R, x)$ defined in \eqref{nu}. The above equality can be rewritten as
	\begin{equation}\label{paraequa2.58}
	\begin{split}
	0 = & \frac{d}{dt}\int_M Y\left(\frac{1}{2}|df|^2+W\right)\varphi^2\, d\textrm{vol}\\
	& +\int_M \left(d^*dY-Y_t+\frac{2Y_r}{r}(1-\nu)\right)\left(\frac{1}{2}|df|^2+W\right)\varphi^2\, d\textrm{vol} \\
	& + \int_M Y\left|f_t+\frac{Y_r}{Y}f_r\right|^2\varphi^2\, d\textrm{vol} +\int_M \left(Y_{rr}-\frac{Y_r}{r}-\frac{Y_r^2}{Y}\right)|f_r|^2\varphi^2\, d\textrm{vol} \\
	& + \int_M Y\langle f_r, f_t+\frac{Y_r}{Y}f_r\rangle 2\varphi\varphi'\, d\textrm{vol} - \int_M Y_r\left(\frac{1}{2}|df|^2+W\right)2\varphi\varphi'\, d\textrm{vol}\\
	&+\int_M Y_r\langle A(x)df, df\rangle\varphi^2\, d\textrm{vol}.
	\end{split}
	\end{equation}
	
	To simplify the equation, we choose $Y$ so that
	\begin{equation}
	Y_{rr}-\frac{Y_r}{r}-\frac{Y_r^2}{Y}=0.
	\end{equation}
	
	One family of solutions is given by
	\begin{equation}
	Y(r, t)=g(t)\exp\left(h(t)r^2\right),
	\end{equation}
	for any positive function $g(t)$ and arbitrary function $h(t)$. 
	
	With this $Y$, we now simplify the expression $d^*dY-Y_t+\dsst \frac{2Y_r}{r}(1-\nu(r,x))$ as much as possible. 
	
	Notice $\dsst d^*dY=-Y_{rr}-\frac{n-1}{r}Y_r+B(r)$ with $B(r)=O(r^2)Y$. Since $B(r)$ will lead to a lower order term, we focus on eliminating $\dsst -Y_{rr}-\frac{m-1}{r}Y_r+\frac{2Y_r}{r}(1-\nu(r,x))$. 
	
	A detailed computation shows
	\begin{equation}
	-Y_{rr}-\frac{m-1}{r}Y_r+\frac{2Y_r}{r}(1-\nu(r,x))=-4h^2r^2Y-(m-2+2\nu(r,x))2hY.
	\end{equation}
	
	Since $m>2$ and $\nu(r,x)$ is small, we require $h(t)\leq 0$ so that we may eliminate $-Y_{rr}-\dsst \frac{m-1}{r}Y_r+\frac{2Y_r}{r}(1-\nu(r,x))$. Therefore, we pick 
	\begin{equation}
	Y(r,t) = g(t)\exp\left(h(t)r^2\right)
	\end{equation}
	with $g(t)> 0$ and $h(t)\leq 0$.
	
	We remove the dependence of the $r$ and $x$ in $\nu(r,x)$ to simplify our computations. By the definition of $\nu_0$, we then get
	\begin{equation}
	\begin{split}
    &	\int_M \left(-Y_{rr}-\frac{m-1}{r}Y_r+\frac{2Y_r}{r}(1-\nu(r,x))\right)\left(\frac{1}{2}|df|^2+W\right)\varphi^2\, d\textrm{vol}\\
     \geq & \int_M \left(-Y_{rr}-\frac{m-1}{r}Y_r+\frac{2Y_r}{r}(1-\nu_0)\right)\left(\frac{1}{2}|df|^2+W\right)\varphi^2\, d\textrm{vol}.
     \end{split}
	\end{equation}
	
	\eqref{paraequa2.58} becomes inequality ($Y_r\leq 0$):
	\begin{equation}\label{parainequal2.64}
	\begin{split}
	0 \geq & \frac{d}{dt}\int_M Y\left(\frac{1}{2}|df|^2+W\right)\varphi^2\, d\textrm{vol}\\
	& +\int_M \left(-Y_t-Y_{rr}-\frac{m-1}{r}Y_r+\frac{2Y_r}{r}(1-\nu_0)\right)\left(\frac{1}{2}|df|^2+W\right)\varphi^2\, d\textrm{vol}\\
	& + \int_M Y\left|f_t+\frac{Y_r}{Y}f_r\right|^2\varphi^2\, d\textrm{vol}  +\int_M B(r)W\varphi^2\, d\textrm{vol}+ \int_M Y_r\langle A(x)df, df\rangle\varphi^2\, d\textrm{vol}\\
	& +\int_M Y\langle f_r, f_t+\frac{Y_r}{Y}f_r\rangle 2\varphi\varphi'\, d\textrm{vol} - \int_M Y_r\left(\frac{1}{2}|df|^2+W\right)2\varphi\varphi'\, d\textrm{vol}.
	\end{split}
	\end{equation}
	
	Imposing the equation on $Y$
	\begin{equation}
	-Y_t-Y_{rr}-\frac{m-1}{r}Y_r+\frac{2Y_r}{r}(1-\nu_0) = 0,
	\end{equation}
	we get
	\begin{equation}
	\left(-4h^2(t)r^2-h'(t)r^2-2h-(m-3+2\nu_0)2h-\frac{g'(t)}{g(t)}\right)Y = 0.
	\end{equation}
	Here $Y>0$ since $g(t)>0$. $h$ and $g$ only depend on $t$. We conclude that
	\begin{equation}
	4h^2(t)+h'(t) = 0,\qquad (m-2+2\nu_0)2h(t)+\frac{g'(t)}{g(t)} = 0.
	\end{equation}
	
	From the first equation, we get
	\begin{equation}
	h(t) = \frac{1}{4(t-C)},
	\end{equation}
	for any constant $C>t$. We then solve for $g(t)$ by plugging this $h$ into the second equation:
	\begin{equation}
	g(t) = D\left|t-C\right|^{1-\frac{m}{2}-\nu_0},
	\end{equation}
	for any $D> 0$. 
	
	Here we pick $D=1$, $C=t_0>t$, then we get
	\begin{equation}
	Y(x,t)=\left(t_0-t\right)^{1-\frac{m}{2}-\nu_0}\exp\left(-\frac{r^2}{4(t_0-t)}\right).
	\end{equation}
	
	We can simplify \eqref{parainequal2.64} into
	\begin{equation}
	\begin{split}
	0 \geq & \frac{d}{dt}\int_{M} Y\left(\frac{1}{2}|df|^2+W\right)\varphi^2\, d\textrm{vol} + \int_M Y\left|f_t+\frac{Y_r}{Y}f_r\right|^2\varphi^2\, d\textrm{vol}\\
	& +\int_M Y_r\langle A(x)df, df\rangle\varphi^2\, d\textrm{vol}+\int_M B(r)W\varphi^2\, dx  +\int_M Y\langle f_r, f_t+\frac{Y_r}{Y}f_r\rangle 2\varphi\varphi'\, d\textrm{vol} \\
	& - \int_M Y_r\left(\frac{1}{2}|df|^2+W\right)2\varphi\varphi'\, d\textrm{vol}\\
	= & \frac{d}{dt}\int_M Y\left(\frac{1}{2}|df|^2+W\right)\varphi^2\, d\textrm{vol} + \int_M Y\left|f_t+\frac{Y_r}{Y}f_r\right|^2\varphi^2\, d\textrm{vol}\\
	& + \text{I} + \text{II} + \text{III} + \text{IV}.
	\end{split}
	\end{equation}
	
	Since $R^2\leq (t_0-t)/4$, we can easily show that
	\begin{equation}
	\left|\textrm{I}\right| = \left|\int_M \frac{r^2}{2(t_0-t)}Y\left(\frac{1}{2}|df|^2+W\right)\varphi^2\, d\textrm{vol}\right|\leq c\int_M Y\left(\frac{1}{2}|df|^2+W\right)\varphi^2\, d\textrm{vol}, 
	\end{equation}
	and
	\begin{equation}
	\left|\text{II}\right| \leq C\int_M \left(\frac{1}{2}|df|^2+W\right)\, d\textrm{vol}\leq CE(f_0).
	\end{equation}
	
	For $\text{III}$, since $R^2\leq t_0-t\leq 4R^2$ with $R<R_M$, we have
	\begin{equation}
	\begin{split}
	|\text{III}| \leq & \frac{1}{2}\int_M Y\left|f_t+\frac{Y_r}{Y}f_r\right|^2\varphi^2\, d\textrm{vol} + C\int_M Y|f_r|^2|\varphi'|^2\,d\textrm{vol} \\
	\leq & \frac{1}{2}\int_M Y\left|f_t+\frac{Y_r}{Y}f_r\right|^2\varphi^2\, d\textrm{vol}+c\int_M Y|f_r|^2\frac{1}{R_M^2}\chi_{[R_M/2, R_M]}\, d\textrm{vol}\\
	\leq & \frac{1}{2}\int_M Y\left|f_t+\frac{Y_r}{Y}f_r\right|^2\varphi^2\, d\textrm{vol} + c\int_M Ye\varphi^2\, d\textrm{vol} + CE(f_0(x)),
	\end{split}
	\end{equation}
	where $\chi_{[R_M/2, R_M]}$ is the indicator of the region $B_{R_M}\setminus B_{R_M/2}$.
	
	\begin{equation}
	\begin{split}
	|\text{IV}| \leq & C\int_M |r\varphi'|\varphi (t_0-t)^{-m/2-\nu_0}\exp\left(-r^2/4(t_0-t)\right)\left(\frac{1}{2}|df|^2+W\right)\, d\textrm{vol}\\
	 \leq & c\int_M Ye\varphi^2\, d\textrm{vol} + CE(f_0(x)).
	\end{split}
	\end{equation}
	
	Here $r$ is only defined when $r$ is close to the boundary of $M$ since that's the only place where $\varphi'$ is supported.
	
	Hence
	\begin{equation}\label{paramonoinequa}
	\begin{split}
	0 \geq & \frac{d}{dt}\int_M Y\left(\frac{1}{2}|df|^2+W\right)\varphi^2\, d\textrm{vol} + \int_M Y\left|f_t+\frac{Y_r}{Y}f_r\right|^2\varphi^2\, d\textrm{vol} \\
	& -\frac{1}{2}\int_M Y\left|f_t+\frac{Y_r}{Y}f_r\right|\varphi^2\,d\textrm{vol}-CE(f_0(x))\\
	\geq & \frac{d}{dt}\int_M Y\left(\frac{1}{2}|df|^2+W\right)\varphi^2\, d\textrm{vol} -c\int_M Y\left(\frac{1}{2}|df|^2+W\right)\varphi^2\, d\textrm{vol} -CE(f_0(x)).
	\end{split}
	\end{equation}
	
	Denote 
	\begin{equation}
	\phi(t) = \int_M Y\left(\frac{1}{2}|df|^2+W\right)\varphi^2\, d\textrm{vol}=\int_M (t_0-t)^{1-\nu_0}G_{z_0}(x,t)\left(\frac{1}{2}|df|^2+W\right)\varphi^2\, d\textrm{vol}.
	\end{equation}
	
	Integrate \eqref{paramonoinequa} in the $t$ direction between $[t_1, t_2]$, where $0<t_1<t_2<t_0$, we have
	\begin{equation}
	\begin{split}
	\phi(t_2)\leq & \exp\left(c(t_2-t_1)\right)\phi(t_1)+\tilde{C}E(f(0,x))\left(\exp\left(c(t_2-t_1)\right)-1\right)\\
	 = & \exp\left(c(t_2-t_1)\right)\phi(t_1)+\hat{C}E(f(0,x))(t_2-t_1).
	 \end{split}
	\end{equation}
	
	Substituting $t_2=t_0-R^2$, and $t_1=t_0-R_0^2$ gives \eqref{phiinequa}.\\
	
	Define
	\begin{equation}
	\Psi(R)= \int_{t_0-4R^2}^{t_0-R^2}\frac{1}{(t_0-t)}\phi(t)\, dt = \int_{t_0-4R^2}^{t_0-R^2}\int_M \frac{1}{(t_0-t)^{\nu_0}} G_{z_0}\left(\frac{1}{2}|df|^2+W\right)\varphi^2\, d\textrm{vol}\, dt.
	\end{equation}
	We have
	\begin{equation}
	\begin{split}
	\frac{d}{dR}\Psi(R) = & \frac{-2R}{R^2}\phi(t_0-R^2)-\frac{-8R}{4R^2}\phi(t_0-4R^2)=\frac{2}{R}\left(\phi(t_0-4R^2)-\phi(t_0-R^2)\right)\\
	& =\frac{2}{R}\left(\phi(t_1)-\phi(t_2)\right).
	\end{split}
	\end{equation}
	
	Now we turn to prove \eqref{psiinequa}. We have
	\begin{equation}
	\phi(t_2)-\phi(t_1)=\int_{t_1}^{t_2}\phi'(t)\, dt\leq \int_{t_1}^{t_2}c\phi(t)+CE_0\, dt =CE_0(t_2-t_1)+c\int_{t_1}^{t_2}\frac{\phi(t)}{t_0-t}(t_0-t)\, dt.
	\end{equation}
	
	So
	\begin{equation}
	\begin{split}
	\frac{d}{dR}\Psi(R) \geq & -\frac{2CE_0}{R}(t_2-t_1)-\frac{c}{R}\int_{t_1}^{t_2}\frac{\phi(t)}{t_0-t}(t_0-t)\, dt\\
	\geq & -\tilde{C}E_0R - \tilde{c}R\Psi(R). 
	\end{split}
	\end{equation}
	
	Integrating from $R$ to $R_0$ (assuming $0< R<R_0<R_M$), we have
	\begin{equation}
	\begin{split}
	\Psi(R)\leq & \Psi(R_0)e^{\frac{\tilde{c}}{2}\left(R_0^2-R^2\right)}+\hat{C}E(f(0, x))\left(e^{\frac{\tilde{c}}{2}(R_0^2-R^2)}-1\right)\\
	\leq & \exp(c(R_0^2-R^2))\Psi(R_0)+\hat{C}E(f(0,x))(R_0^2-R^2).
	\end{split}
	\end{equation}
	
\end{proof}


%% file: Chapter3/Chapter3.tex
\def\dsst{\displaystyle}
\def\R{\mathds{R}}
\chapter{Moser Iteration}
\label{chap:moser}

In this chapter we recall the Moser iteration arguments that allow us to derive $L^\infty$ bounds for smooth functions $e$ satisfying either an elliptic inequality
$$\Delta e \leq Ce,$$
or a parabolic inequality
$$(\frac{d}{dt}+\Delta) e \leq Ce.$$

The contents of this section are not original, but will be used in the next section to obtain cleaner and more effective $\epsilon$-regularity arguments than we have found in the literature (see \cite{schoen1984analytic} and \cite{StruweSigma}).

\section{Moser Iteration for Elliptic subsolutions}



\begin{lemma}
	
	Let $M$ be a compact manifold. Let $e$ be a positive function from $M\rightarrow \R^1$ satisfying
	\begin{equation}\label{Ellipticinequ}
	\Delta e\leq C_0e,
	\end{equation}  
	where $\Delta$ is the Laplace-Beltrami operator on $M$. 
	
	Then
	\begin{equation}\label{ellipticMoser}
	 \left\|e\right\|_{L^\infty(B_{\delta R}(x_0))}\leq \widetilde{C}(M)\left(C_0+\frac{16}{(R-\delta R)^2}\right)^{\frac{m}{2}}\int_{B_R(x_0)} e\, d\textup{vol},
	 \end{equation}
	for $0\leq\delta<1$ and any geodesic ball $B_R(x_0)\in M$ where $R<R_M$. Here $\widetilde{C}(M)$ is a constant only depends on $M$.
	
	Particularly, when $\delta=0$, \eqref{ellipticMoser} implies a single point estimate of $e(f)$:
	\begin{equation}
	e(f)(x_0)\leq \tilde{C}(M)\left(C_0+\frac{16}{R^2}\right)^{\frac{m}{2}}\int_{B_{R}(x_0)}e\,d\textup{vol}.
	\end{equation}
	
	If $\dsst \sup_M e(f)$ exists and is achieved at $x_0$, then 
	\begin{equation}
	\sup_M e(f) \leq \tilde{C}(M)\left(C_0+\frac{16}{R^2}\right)^{\frac{m}{2}}\int_{B_R(x_0)} e\,d\textup{vol},
	\end{equation}
	for any $R<R_M$. 
\end{lemma}

\begin{proof}
	
	Without loss of generality, we may assume $x_0=0$ and denote $B_R(0)=B_R$. We define the radial cut-off function $\eta(s)$ such that $\eta(s)=1$ for $s\leq \frac{3}{2}$, and $\eta(s)=0$ for $s\geq 2$, with $|\eta'(s)|\leq 4$. Define radial function $\eta_k(x)=\eta_k(|x|)=\eta\left(\dsst 1+2^{k-1}\left(\frac{|x|}{R}-1\right)\right)$. We can show that $\eta_k(x)=1$ in $B_{r_k}$ where $r_{k}=\dsst R\left(1+\frac{1}{2^{k}}\right)$, and $\eta(x)$ is supported in $B_{r_{k-1}}$, $\dsst \left|d\eta_k\right|\leq \frac{2^{k+1}}{R}$, and $d\eta_k$ is supported in the set $B_{r_{k-1}}\setminus B_{r_{k}}$. \\

	Multiply \eqref{Ellipticinequ} by $\eta_k e^{p_k-1}$ on both sides, and integrate by parts to get
	\begin{equation}\label{EllipticMoser2}
	\begin{split}
	C_0\int_{B_R}|\eta_ke^{p_k}|^2 \geq & \langle \Delta e, \eta_k^2e^{2p_k-1}\rangle_{L^2}=\langle de, d(\eta_k^2e^{2p_k-1})\rangle_{L^2}\\
	= &2\langle (\eta_ke^{p_k-1})de, ed(\eta_ke^{p_k-1})\rangle_{L^2} +\int_{B_R}\eta_k^2e^{2p_k-2}|de|^2\\
	= & \int_{B_R} \langle d(\eta_ke^{p_k}), d(\eta_ke^{p_k})\rangle d\textrm{vol} -\int_{B_R}\langle d(\eta_ke^{p_k-1})e, d(\eta_ke^{p_k-1})e\rangle d\textrm{vol}\\
	= & \left\|d(\eta_ke^{p_k})\right\|^2_{L^2(B_R)}-\left\|d(\eta_ke^{p_k-1})e\right\|_{L^2(B_R)}^2.
	\end{split}
	\end{equation}	
	
	Therefore, 
	\begin{equation}
	\begin{split}
	\left\|d(\eta_k e^{p_k})\right\|^2_{L^2(B_R)}\leq &  C_0\left\|\eta_ke^{p_k}\right\|_{L^2(B_R)}^2+\left\|d(\eta_ke^{p_k-1})e\right\|_{L^2(B_R)}^2\\
	= & C_0\left\|\eta_ke^{p_k}\right\|_{L^2(B_R)}^2+\left\|e^{p_k}d\eta_k+(p_k-1)\eta_ke^{p_k-1}de\right\|_{L^2(B_R)}^2\\
	= & C_0\left\|\eta_ke^{p_k}\right\|_{L^2(B_R)}^2+\left\|\frac{p_k-1}{p_k}\eta_kd(e^{p_k})+e^{p_k}d\eta_k\right\|_{L^2(B_R)}^2\\
	= & C_0\left\|\eta_ke^{p_k}\right\|_{L^2(B_R)}^2+\left\|\frac{p_k-1}{p_k}d(\eta_ke^{p_k})+\frac{1}{p_k}e^{p_k}d\eta_k\right\|_{L^2(B_R)}^2\\ 
	= & C_0\left\|\eta_ke^{p_k}\right\|_{L^2(B_R)}^2+\left(\frac{p_k-1}{p_k}\right)^2\left\|d(\eta_ke^{p_k})\right\|^2_{L^2(B_R)}\\
	& +2\frac{p_k-1}{p_k^2}\langle d(\eta_ke^{p_k}),e^{p_k}d\eta_k\rangle_{L^2(B_R)}+\frac{1}{p_k^2}\left\|e^{p_k}d\eta_k\right\|^2_{L^2(B_R)}.
	\end{split}
	\end{equation}
	
	We apply Cauchy-Schwarz to the cross term to get,
	\begin{equation}
	\begin{split}
	& \frac{2p_k-1}{p_k^2}\left\|d(\eta_ke^{p_k})\right\|_{L^2(B_R)}^2 \\
	\leq &   C_0\left\|\eta_ke^{p_k}\right\|^2_{L^2(B_R)}+\frac{p_k-1}{p_k^2}\left\|d(\eta_ke^{p_k})\right\|^2_{L^2(B_R)}+\frac{1}{p_k^2}\left\|e^{p_k}d\eta_k\right\|^2_{L^2(B_R)}.
	\end{split}
	\end{equation}
	
	Therefore,
	\begin{equation}
	\frac{1}{p_k}\left\|d(\eta_ke^{p_k})\right\|^2_{L^2(B_R)} \leq C_0\left\|\eta_ke^{p_k}\right\|_{L^2(B_R)}^2+\frac{1}{p_k^2}\left\|e^{p_k}d\eta_k\right\|^2_{L^2(B_R)}.
	\end{equation}
	
	That is,
	\begin{equation}
	\begin{split}
	\left\|d(\eta_ke^{p_k})\right\|_{L^2(B_R)}^2\leq & C_0p_k\left\|\eta_ke^{p_k}\right\|^2_{L^2(B_R)}+\left\|e^{p_k}d\eta_k\right\|^2_{L^2(B_R)}\\
	\leq & C_0p_k\left\|\eta_ke^{p_k}\right\|_{L^2(B_R)}^2+\frac{4^{k+1}}{r^2}\left\|e^{p_k}\chi_{k-1}\right\|_{L^2(B_R)}^2,
	\end{split}
	\end{equation}
	where $\chi_{k-1}$ is the indicator function of $ B_{r_{k-1}}\setminus B_{r_k}$.
	
	Since $M$ is compact, from the Sobolev inequality \cite[Theorem 2.6]{hebey2000nonlinear} we have
	\begin{equation}
	\left\|d(\eta_k e^{p_k})\right\|_{L^2(B_R)}\geq C(M) \left\|\eta_k e^{p_k}\right\|_{L^{\frac{2m}{m-2}}(B_R)},
	\end{equation}
	where $C(M)$ is a constant depending only $M$.\\
	
	Setting $p_1 = 1$ and $p_{k+1}=p_k\frac{m}{m-2}$, we have
	\begin{equation}
	\begin{split}
	C_0p_k\left\|e\right\|_{L^{2p_k}(B_{r_{k-1}})}^{2p_k} = & C_0p_k\int_{B_{r_{k-1}}} e^{2p_k}\, d\textrm{vol}\geq C_0p_k\left\|\eta_k e^{p_k}\right\|_{L^2(B_R)}^2 \\
	\geq & \left\|d(\eta_ke^{p_k})\right\|_{L^2}^2-\frac{2^{k+2}}{r^2}\left\|e^{p_k}\chi_{k-1}\right\|_{L^2(B_R)}^2\\
	\geq & C(M)\left\|\eta_ke^{p_k}\right\|^2_{L^{\frac{2m}{m-2}}}-\frac{4^{k+1}}{r^2}\left\|e^{p_k}\chi_{k-1}\right\|_{L^2(B_R)}^2\\
	= & C(M)\left(\int_{B_{r_k}} e^{2p_k\frac{m}{m-2}}d\textrm{vol}\right)^{\frac{m-2}{m}}-\frac{4^{k+1}}{r^2}\left\|e^{p_k}\chi_{k-1}\right\|_{L^2(B_R)}^2\\
	= & C(M)\left(\int_{B_{r_k}} e^{2p_{k+1}}d\textrm{vol}\right)^{\frac{2p_k}{2p_{k+1}}}-\frac{4^{k+1}}{r^2}\left\|e^{p_k}\chi_{k-1}\right\|^2_{L^2(B_R)}.
	\end{split}
	\end{equation}
	
	Rearranging, we obtain
	\begin{equation}
	\left(C_0p_k+4^k\frac{4}{r^2}\right)\left\|e\right\|_{L^{2p_k}(B_{r_{k-1}})}^{2p_k}\geq C(M)\left\|e\right\|_{L^{2p_{k+1}}(B_{r_k})}^{2p_k}.
	\end{equation}
	
	That is,
	\begin{equation}
	\left(\frac{C_0p_k}{C(M)}+4^k\frac{4}{C(M)r^2}\right)^{\frac{1}{p_k}}\left\|e\right\|^2_{L^{2p_k}(B_{r_{k-1}})}\geq \left\|e\right\|^2_{L^{2p_{k+1}}(B_{r_k})}.
	\end{equation}
	
	Iterating, we get
	\begin{equation}\label{iterate}
	\prod_{k=1}^\infty \left(\frac{C_0p_k}{C(M)}+4^k\frac{4}{C(M)r^2}\right)^{\frac{1}{p_k}}\left\|e\right\|_{L^2(B_{2R})}^2\geq \left\|e\right\|_{L^\infty(B_R)}^2.
	\end{equation}
	
	
	Since $\dsst p_k=\left(\frac{m}{m-2}\right)^{k-1}$, $p_k\leq 4^{k-1}$. The product factor on the left-hand side of \eqref{iterate} can be bounded by
	\begin{equation}
	\begin{split}
	\prod_{k=1}^\infty \left(\frac{C_0p_k}{C(M)}+4^{k-1}\frac{16}{r^2}\right)\leq & \prod_{k=1}^\infty \left[\frac{4^k}{C(M)}\left(C_0+\frac{16}{r^2}\right)\right]^{\frac{1}{p_k}}\\
	= & 4^{\frac{m^2}{4}}\left(\frac{1}{C(M)}\right)^{\frac{m}{2}}\left(C_0+\frac{16}{r^2}\right)^{\frac{m}{2}}\\
	= &\widetilde{C}(M)\left(C_0+\frac{16}{r^2}\right)^{\frac{m}{2}}.
	\end{split}
	\end{equation}
	
	Hence,
	\begin{equation}\label{2RtoR}
	\left\|e\right\|_{L^\infty(B_{R})}^2\leq \tilde{C}(M)\left(C_0+\frac{16}{R^2}\right)^{\frac{m}{2}} \int_{B_{2R}} e^2\, d\textrm{vol}.
	\end{equation}
	
	
	To pass from a bound in terms of $\|e\|_{L^2(B_{2R})}$ to one in terms of $\|e\|_{L^2(B_R)}$, we first consider  $r_k=R\left(\delta + \frac{1-\delta}{2^{k}}\right)$. This changes the bound of the derivative of $\eta_k$ to $|d\eta_k|\leq \dsst \frac{2\cdot 2^k}{R(1-\delta)}$. Then \eqref{2RtoR} will turn into
	\begin{equation}\label{deltaRtoR}
	\sup_{x\in B_{\delta R}(x_0)}e^2\leq \tilde{C}(M)\left(C_0+\frac{16}{(R-\delta R)^2}\right)^{\frac{m}{2}}\int_{B_R(x_0)}e^2\,d\textrm{vol}.
	\end{equation}
	
	To show that $\sup e$ can be bounded by its $L^1$ norm, we use the following lemma from \cite[Lemma 4.1]{lovelylemma}
	\begin{lemma}\label{speciallemma}
		Let $\varphi(t)$ be a bounded nonnegative function defined on the interval $[T_0, T_1]$, where $0\leq T_0<T_1$. Suppose that for any $T_0\leq t<s\leq T_1$, $\varphi$ satisfies
		\begin{equation}
		\varphi(t)\leq \theta\varphi(s)+\frac{A}{(s-t)^{\alpha}} + B,
		\end{equation}
		where $\theta, A, B$ and $\alpha$ are nonnegative constants, $\theta<1$. Then
		\begin{equation}
		\varphi(\rho)\leq C\left[\frac{A}{(R-\rho)^{\alpha}}+B\right],\qquad \forall T_0\leq \rho<R\leq T_1,
		\end{equation}
		where $C$ depends only on $\alpha, \theta$. 
	\end{lemma}
	
	We apply this lemma by setting $t=\delta R$, $s=R$, and $\varphi(t)=\dsst \sup_{B_{t}(x_0)}e$. By Cauchy-Schwarz and \eqref{deltaRtoR}, we get
	\begin{equation}
	\begin{split}
	\varphi(t) = & \sup_{B_{\delta R}(x_0)}e\\
	\leq & \left(C_0\tilde{C}+\frac{16\tilde{C}}{(R-\delta R)^2}\right)^{\frac{m}{4}}\sqrt{\int_{B_R(x_0)}e^2\, d\textrm{vol}}\\
	\leq & \sqrt{\sup_{B_R(x_0)}e}\left(C_0\tilde{C}+\frac{16\tilde{C}}{(R-\delta R)^2}\right)^{\frac{m}{4}}\sqrt{\int_{B_R(x_0)}e\, d\textrm{vol}}\\
	\leq & \frac{1}{2}\sup_{B_R(x_0)}e+2\left(\tilde{C}_0+\frac{16\tilde{C}}{(R-\delta R)^2}\right)^{\frac{m}{2}}\int_{B_R(x_0)}e\, d\textrm{vol}\\
	= & \frac{1}{2}\varphi(s)+\frac{A}{(s-t)^m}\int_{B_R(x_0)}e\, d\textrm{vol},
	\end{split}
	\end{equation}
	where $A$ can be taken as a constant slightly larger than $\dsst \tilde{C}_0(R-\delta R)^2+16\tilde{C}$ since $(R-\delta R)^2$ is bounded. 
	
	Applying the lemma, we have
	
	\begin{equation}
	\varphi(t)\leq C\left[\frac{A}{(s-t)^m}\int_{B_R(x_0)} e\, d\textrm{vol}\right].
	\end{equation}

	That is,
	\begin{equation}
	\sup_{B_{\delta R}(x_0)} e\leq \left(C_1(M)C_0+\frac{C_2(M)}{(R-\delta R)^2}\right)^{\frac{m}{2}}\int_{B_R(x_0)}e\, d\textrm{vol}.
	\end{equation}
	
	In particular, we have the single point estimate at the center $x_0$ (a special case when $\delta =0$):
	\begin{equation}
	e(f)(x_0)\leq \left(C_1(M)C_0 + \frac{C_2(M)}{R^2}\right)^{\frac{m}{2}}\int_{B_{R}(x_0)} e(f)\, d\textrm{vol}.
	\end{equation}
	
	If $\dsst \sup_M e(f)=e(x_1)$, then from the Moser iteration, we can easily get the result
	\begin{equation}
	\sup_M e(f) \leq \left(C_1(M)+\frac{C_2(M)}{R^2}\right)^{\frac{m}{2}} \int_{B_{R}(x)} e(f)\, d\textrm{vol},
	\end{equation}
	for $x\in M$ such that $e(f)(x)=\dsst \sup_M e(f)$ and any $R<R_M$ with $B_R(x)\subset M$.
	
\end{proof}

\section{Moser-Harnack Iteration for Parabolic Subsolutions}

\begin{lemma}
	
	Let $z_0=(x_0, t_0)\in M\times I$ with $M$ compact. Suppose a positive function $e:M\rightarrow \R^1$ satisfies
	\begin{equation}\label{parasubsolution}
	\left(\frac{\partial}{\partial t}+\Delta\right)e\leq C_0 e \qquad \text{in }P_R(z_0),
	\end{equation}
	for some $C_0\geq 0$ where it may depend on $e$, and $P_R(z_0)$ is defined to be the cylinder \begin{equation}
	P_R(z_0) = \left\{(x, t)\in M \times I\,| |x-x_0|\leq R, |t-t_0|\leq R^2\right\}
	\end{equation}
	for $R<R_M$.\\
	
	Then
	\begin{equation}
	\sup_{P_{\delta R}(z_0)}e \leq  \left(C_1(M)C_0+\frac{C_2(M)}{(R-\delta R)^2 }\right)^{\frac{m+2}{2}}\int_{P_R(z_0)} e\, d\textup{vol}\,dt,
	\end{equation}
	where $C_1(M)$ and $C_2(M)$ only depend on $M$.\\
	
	In particular, when $\delta=0$, we have
	\begin{equation}\label{singlepoint}
	e(f)(z_0) \leq \left(C_1(M)C_0+\frac{C_2(M)}{ R^{2}}\right)^{\frac{m+2}{2}}\int_{P_R(z_0)} e\, d\textup{vol}\,dt.
	\end{equation}
	
	If $\dsst \sup_{M\times I} e(f)$ exists, then 
	\begin{equation}\label{entireM}
	\sup_{M\times I} e(f) \leq \left(C_1(M)C_0+\frac{C_2(M)}{ R^{2}}\right)^{\frac{m+2}{2}}\int_{P_R(z_0)} e\, d\textup{vol}\,dt,
	\end{equation}
	for any $R<R_M$.
\end{lemma}

\begin{proof}
	
	Without loss of generality, we may assume $z_0=(x_0, t_0)=(0, 0)$. We denote $P_R(0)=P_R$ and so on. First we assume
	\begin{equation}
	\left(\frac{\partial}{\partial t}+\Delta\right)e(f) \leq C_0e(f) \quad \text{in }P_{2R}(z_0).
	\end{equation}
	
	Define the spatial cut-off functions $\eta_k$ as in the elliptic case. Similarly, we can also define cut-off functions in the time direction. Let $g_k(t)=1$ for $|t|\leq R^2\left(\dsst 1+\frac{1}{2^k}\right)$ and $g_k(t)=0$ for $|t|\geq R^2\left(\dsst 1+\frac{2}{2^{k-1}}\right)$ so that $|g_k'(t)|\leq \dsst \frac{2^k}{R^2}$. \\
	
	Multiply the equation by $e^{2p_n-1}\eta_n^2g_n^2$ on both sides, then integrate over $P_{2R}$. With integration by parts, we have
	\begin{equation}
	\begin{split}
	 & C_0\int_{P_{2R}} e^{2p_n}\eta_n^2g_n^2\, d\textrm{vol}\,dt\\
	  \geq & \int_{P_{2R}}\langle \frac{\partial e}{\partial t}, e^{2p_n-1}\eta_n^2g_n^2\rangle\, d\textrm{vol}\,dt+\int_{P_{2R}}\langle \Delta e, e^{2p_n-1}\eta_n^2g_n^2\rangle\, d\textrm{vol}\, dt\\
	= & \int_{P_{2R}}\langle \frac{\partial e}{\partial t}, e^{2p_n-1}\eta_n^2g_n^2\rangle \, d\textrm{vol}\,dt+\int_{P_{2R}}\langle de, d\left(e^{2p_n-1}\eta_n^2\right)g_n^2\rangle\, d\textrm{vol}\, dt \\
	= & \frac{1}{2p_n} \int_{P_{2R}} \frac{d}{dt}\left(e^{2p_n}\right)\eta_n^2g_n^2\, d\textrm{vol}\, dt+\frac{2p_n-1}{p_n^2}\int_{P_{2R}}|d(e^{p_n})|^2\eta_n^2g_n^2\, d\textrm{vol}\, dt\\
	& +\frac{1}{p_n}\int_{P_{2R}} e^{p_n}g_n^2\langle d(e^{p_n}), d(\eta_n^2)\rangle\, d\textrm{vol}\, dt\\
	= & \frac{1}{2p_n}\int_{P_{2R}}\frac{d}{dt}\left(\left(e^{2p_n}\eta_n^2g_n^2\right)\right)\, d\textrm{vol}\, dt-\frac{1}{2p_n}\int_{P_{2R}}e^{2p_n}\eta_n^2\frac{d}{dt}\left(g_n^2\right)\, d\textrm{vol}\, dt\\
	& + \frac{2p_n-1}{p_n^2}\int_{P_{2R}}|d(e^{p_n})|^2\eta_n^2g_n^2\, d\textrm{vol}\, dt+\frac{1}{p_n}\int_{P_{2R}}e^{p_n}g_n^2\langle d(e^{p_n}), d(\eta_n^2)\rangle\, d\textrm{vol}\, dt.
	\end{split}
	\end{equation}
	
	
	We bound the last two terms by $\left\|d(e^{p_n}\eta_ng_n)\right\|_{L^2(P_{2R})}$, and pick $p_1=1,\ \dsst p_k=\left(\frac{m+2}{m}\right)p_{k-1}\geq 1$
	
	\begin{equation}
	\begin{split}
	& C_0\int_{P_{2R}} e^{2p_n}\eta_n^2g_n^2\, d\textrm{vol}\,dt\\
	\geq  & \frac{1}{2p_n}\int_{P_{2R}}\frac{d}{dt}\left(e^{2p_n}\eta_n^2g_n^2\right)\, d\textrm{vol}\, dt-\frac{1}{2p_n}\int_{P_{2R}}e^{2p_n}\eta_n^2\frac{d}{dt}\left(g_n^2\right)\, d\textrm{vol}\, dt\\
	& +\frac{1}{p_n}\int_{P_{2R}}\left|d\left(e^{p_n}\eta_ng_n\right)\right|^2\, d\textrm{vol}\, dt+\frac{p_n-1}{p_n^2}\int_{P_{2R}}|d(e^{p_n})|^2\eta_n^2g_n^2\, d\textrm{vol}\, dt\\
	& -\frac{1}{p_n}\int_{P_{2R}}f^{2p_n}g_n^2\left|d(\eta_n)\right|^2\, d\textrm{vol}\, dt\\
	\geq & \frac{1}{2p_n}\int_{P_{2R}}\frac{d}{dt}\left(e^{2p_n}\eta_n^2g_n^2\right)\, d\textrm{vol}\, dt +  \frac{1}{p_n}\int_{P_{2R}}\left|d\left(e^{p_n}\eta_ng_n\right)\right|^2\, d\textrm{vol}\, dt\\
	& -\frac{1}{p_n}\int_{P_{2R}} e^{2p_n}\left(\frac{1}{2}\eta_n^2\frac{d}{dt}\left(g_n^2\right)+g_n^2\left|d\left(\eta_n\right)\right|^2\right)\, d\textrm{vol}\, dt.
	\end{split}
	\end{equation}
	
	Rearranging, we have
	\begin{equation}
	\begin{split}
	 & \int_{P_{2R}}\frac{d}{dt}\left(e^{2p_n}\eta_n^2g_n^2\right)\, d\textrm{vol}\, dt+ 2\int_{P_{2R}}\left|d\left(e^{p_n}\eta_ng_n\right)\right|^2\, d\textrm{vol}\, dt\\
	 \leq &\int_{P_{2R}} e^{2p_n}\left(2C_0p_n\eta_n^2g_n^2+\eta_n^2\frac{d}{dt}\left(g_n^2\right)+2g_n^2\left|d\left(\eta_n\right)\right|^2\right)\, d\textrm{vol}\, dt.
	\end{split}
	\end{equation}
	
	The preceding inequality still holds if we replace $P_{2R}$ by $P_{2R}^\tau:= B_{2R}\times[-4R^2, \tau) $ for any $\tau\in (-4R^2, 0]$. Since $g_n^2=0$ when $t=-4R^2$, integration by parts in the $t$ direction yields
	\begin{equation}
	\begin{split}
	& \int_{B_{2R}\times \{\tau\}} \eta_n^2g_n^2e^{2p_n}\, d\textrm{vol}+2\int_{P_{2R}}\left|d\left(e^{p_n}\eta_ng_n\right)\right|^2\, d\textrm{vol}\, dt\\
	\leq &\int_{P_{2R}} e^{2p_n}\left(2C_0p_n\eta_n^2g_n^2+\eta_n^2\frac{d}{dt}\left(g_n^2\right)+2g_n^2\left|d\left(\eta_n\right)\right|^2\right)\, d\textrm{vol}\, dt.
	\end{split}
	\end{equation}
	
	By Young's inequality, for any $a,\ b\geq0$, and $\dsst \frac{1}{p}+\frac{1}{q}=1$, 
	\begin{equation}
	\frac{a^p}{p}+\frac{b^q}{q}\geq ab.
	\end{equation}
	Pick $p=\dsst \frac{m+2}{2}$ and $q=\dsst \frac{m+2}{m}$ and view $\dsst \int_{B_{2R}\times \{\tau\}}e^{2p_n}\eta_n^2g_n^2\, d\textrm{vol}\, dt=\frac{a^p}{p}$, \\$\dsst 2\int_{P_{2R}}\left|d\left(e^{p_n}\eta_ng_n\right)\right|^2\, d\textrm{vol}\, dt=\frac{b^q}{q}$, we have
	
	\begin{equation}
	\begin{split}
	& \left(\frac{m+2}{2}\int_{B_{2R}\times\{\tau\}}e^{2p_n}\eta_n^2g_n^2\, d\textrm{vol}\right)^{\frac{1}{p}}\left(\frac{2m}{m+2}\int_{P_{2R}}\left|d\left(e^{p_n}\eta_ng_n\right)\right|^2\, d\textrm{vol}\, dt\right)^{\frac{1}{q}}\\
	\leq & \int_{P_{2R}} e^{2p_n}\left(2C_0p_n\eta_n^2g_n^2+\eta_n^2\frac{d}{dt}\left(g_n^2\right)+2g_n^2\left|d\left(\eta_n\right)\right|^2\right)\, d\textrm{vol}\, dt.
	\end{split}
	\end{equation}
	
	This works for any $\tau\in(-4r^2, 0]$, therefore
	\begin{equation}
	\begin{split}
	& \left(\sup_{\tau\in(-4R^2, 0]}\frac{m+2}{2}\int_{B_{2R}\times\{\tau\}}e^{2p_n}\eta_n^2g_n^2\, dx\right)^{\frac{1}{p}}\left(\frac{2m}{m+2}\int_{P_{2R}}\left|d\left(e^{p_n}\eta_ng_n\right)\right|^2\, d\textrm{vol}\, dt\right)^{\frac{1}{q}}\\
	\leq & \int_{P_{2R}} e^{2p_n}\left(2C_0p_n\eta_n^2g_n^2+\eta_n^2\frac{d}{dt}\left(g_n^2\right)+2g_n^2\left|d\left(\eta_n\right)\right|^2\right)\, d\textrm{vol}\, dt.
	\end{split}
	\end{equation}
	
	By the Sobolev Lemma, for $I$ an interval in $\R^1$, 
	\begin{equation}\label{sobolevu}
	\int_{B_{\rho}\times I} u^{2\frac{m+2}{m}}\, d\textrm{vol}\, dt \leq C(M)\left(\sup_{\tau\in I}\int_{B_\rho\times\{t\}} u^2\, d\textrm{vol}\right)^{\frac{2}{m}}\times \int_{B_\rho\times I}\left|du\right|^2\, d\textrm{vol}\, dt
	\end{equation}
	
	By \eqref{sobolevu} with $u=e^{p_n}\eta_ng_n$, we have
	\begin{equation}
	\begin{split}
	& \widetilde{C}(M)\left(\int_{P_{2R}}\left(e^{p_n}\eta_ng_n\right)^{2\frac{m+2}{m}}\, d\textrm{vol}\, dt\right)^{\frac{m}{m+2}}\\
	\leq & \left[\left(\sup_{\tau\in(-4R^2, 0]}\int_{B_{2R}\times\{\tau\}}(e^{p_n}\eta_ng_n)^2\right)^{\frac{2}{m}}\times \int_{P_{2R}}\left|d(e^{p_n}\eta_ng_n)\right|^2\, d\textrm{vol}\, dt\right]^{\frac{m}{m+2}}\\
	\leq & \int_{P_{2R}} e^{2p_n}\left(2C_0p_n\eta_n^2g_n^2+\eta_n^2\frac{d}{dt}\left(g_n^2\right)+2g_n^2\left|d\left(\eta_n\right)\right|^2\right)\, d\textrm{vol}\, dt.
	\end{split}
	\end{equation}
	
	Since $\dsst |d(g_n)|\leq \frac{2^n}{R^2},\ |d(\eta_n)|\leq \frac{2^{n+1}}{R}$, and $g_n$, $\eta_n$ are cutoff functions both compactly supported by $P_{R\left(1+\frac{1}{2^{n-1}}\right)}$, we can use the information to restrict the domain of the integrals on both sides
	\begin{equation}
	\widetilde{C}(M)\left(\int_{P_{R\left(1+\frac{1}{2^n}\right)}}e^{2p_n\frac{m+2}{m}}\, d\textrm{vol}\, dt\right)^{\frac{m}{m+2}} \leq \left(2C_0p_n+\frac{4^{n+1}}{R^2}\right)\int_{P_{R\left(1+\frac{1}{2^{n-1}}\right)}} e^{2p_n}\, d\textrm{vol}\, dt.
	\end{equation}
	
	That is
	\begin{equation}
	\widetilde{C}(M)\left\|e\right\|_{L^{2p_{n+1}}(P_{r_n})}^{2p_n} \leq \left(C_0p_n+\frac{4^n}{R^2}\right)\left\|e\right\|_{L^{2p_n}(P_{r_{n-1}})}^{2p_n}.
	\end{equation}
	
	Hence,
	\begin{equation}
	\left\|e\right\|^{2}_{L^{2p_{n+1}}(P_{r_n})}\leq \left(\frac{C_0p_n}{\widetilde{C}(M)}+\frac{4^n}{\widetilde{C}(M)R^2}\right)^{\frac{1}{p_n}}\left\|e\right\|^2_{L^{2p_n}(P_{r_{n-1}})}.
	\end{equation}
	
	Iterate until $p_n\rightarrow\infty$, we have
	\begin{equation}
	\begin{split}
	\left\|e\right\|^2_{L^\infty(P_R)}\leq &  \prod_{n=1}^\infty\left(\frac{2^{n}}{\widetilde{C}(M)}\left(C_0+\frac{1}{R^2}\right)\right)^{\frac{1}{p_n}}\left\|e\right\|^2_{L^2(P_{2R})} \\
	= & \left(\frac{4}{\widetilde{C}(M)}\left(C_0+\frac{1}{R^2}\right)\right)^{\frac{m+2}{2}}\left\|e\right\|^2_{L^2(P_{2R})}\\
	= & \left(C_1(M)C_0+\frac{C_2(M)}{R^2}\right)^{\frac{m+2}{2}}\left\|e\right\|_{L^2(P_{2R})}^2,
	\end{split}
	\end{equation}
	for $C_1(M),\ C_2(M)$ being just the Sobolev constants that only depend on $M$.\\
	
	As for the general case when \eqref{parasubsolution} is true in $P_R(z_0)$, we have
	\begin{equation}
	\left(\sup_{z\in P_{\delta R}(z_0)}e\right)^2 \leq \left(C_1(M)C_0+\frac{C_2(M)}{(R-\delta R)^2}\right)^{\frac{m+2}{2}}\int_{P_{R}(z_0)} e^2\, d\textrm{vol}\, dt.
	\end{equation}
	
	We then apply Lemma \ref{speciallemma} in order to get the $L^\infty$ norm bounded by the $L^1$ norm of $e(f)$ in a larger cylinder.
	
	We use the same trick as in the elliptic case. From Cauchy-Schwarz, we have
	\begin{equation}
	\begin{split}
	\sup_{P_{\delta R}(z_0)} e(f) \leq &  \left(C_1(M)C_0+\frac{C_2(M)}{(R-\delta R)^2}\right)^{\frac{m+2}{4}}\sqrt{\int_{P_R(z_0)} e^2\, d\textrm{vol}\, dt}\\
	\leq & \left(C_1(M)C_0+\frac{C_2(M)}{(R-\delta R)^2}\right)^{\frac{m+2}{4}}\sqrt{\sup_{P_{R(z_0)}} e(f)}\sqrt{\int_{P_{R}(z_0)} e(f)\, d\textrm{vol}\, dt} \\
	\leq & \frac{1}{2}\sup_{P_{R(z_0)}} e(f)+\left(C_1(M)C_0+\frac{C_2(M)}{(R-\delta R)^2}\right)^{\frac{m+2}{2}}\int_{P_{R}(z_0)} e(f)\, d\textrm{vol}\, dt.
	\end{split}
	\end{equation}
	Now apply the lemma by setting $t=\delta R$, $s=R$, and 
	$$ \varphi(t)=\sup_{P_t(z_0)} e(f), $$
	we have
	$$ \varphi(t) \leq \frac{1}{2}\varphi(s)+\left(C_1(M)C_0+\frac{C_2(M)}{(R-\delta R)^2}\right)^{\frac{m+2}{2}}\int_{P_{R}(z_0)} e(f)\, d\textrm{vol}\, dt. $$
	
	Since $C_1(M)C_0$ is bounded, we get
	$$ \varphi(t)\leq \frac{1}{2}\varphi(s)+\frac{A}{(s-t)^{m+2}}+B, $$
	where
	$$ A=\tilde{C}_2(M)\int_{P_{R}(z_0)} e(f)\, d\textrm{vol}\, dt,\qquad B=\tilde{C}_1(M)C_0\int_{P_{R}(z_0)} e(f)\, d\textrm{vol}\, dt. $$
	
	From the lemma, we get
	\begin{equation}
	\begin{split}
	\varphi(t) = &\sup_{P_{\delta R}(z_0)} e(f) \\
	\leq & C\left[\frac{A}{(R-\delta R)^{m+2}}+B\right] \\
	\leq & \left(\tilde{C}_1(M)C_0+\frac{\tilde{C}_2(M)}{(R-\delta R)^2}\right)^{\frac{m+2}{2}}\int_{P_{R}(z_0)} e(f)\,d\textrm{vol}\, dt,
	\end{split}
	\end{equation}
	for any $0\leq \delta < 1$. \\
	
	In particular, we have
	\begin{equation}
	e(f)(z_0) \leq \left(\tilde{C}_1(M)C_0+\frac{\tilde{C}_2(M)}{R^2}\right)^{\frac{m+2}{2}}\int_{P_{R}(z_0)} e(f)\, d\textrm{vol}\, dt, 
	\end{equation}
	for any $R$ where \eqref{parasubsolution} holds. \\
	
	Suppose $\dsst \sup_{M\times I} e(f)$ exists and is achieved at $z_0$, we also have
	\begin{equation}
	\sup_{M\times I} e(f) \leq \left(\tilde{C}_1(M)C_0+\frac{\tilde{C}_2(M)}{R^2}\right)^{\frac{m+2}{2}}\int_{P_{R}(z_0)} e(f)\, d\textrm{vol}\, dt,
	\end{equation}
	for any $R$.\\

\end{proof}

%% file: Chapter4/Chapter4.tex
\chapter{$\epsilon$-Regularity}
\label{ch4:epsilon}

In this chapter, we  apply  Moser iteration and the monotonicity formulas derived in Chapter \ref{chap:moser}   to derive  new $\epsilon$-regularity results.  We obtain better bounds for $\sup e$ under suitable hypotheses on  the ratios $\mu_0$, $\nu_0$, and  $p_0$. This proof is derived directly from Moser iteration without introducing any rescaling of the function $e$. Therefore, we can see clearly how $\epsilon_0$ depends on the constants playing a role in the Moser iteration. 

\section{Elliptic $\epsilon$-Regularity}

We obtain estimates for positive functions $e:B_{R}(x_0)\rightarrow \R^1$ which satisfy a monotonicity relation and the elliptic inequality
\begin{equation}\label{ellipinequal}
\Delta e\leq Ce^2.
\end{equation}

\begin{theorem}\label{epsilonelliptic}
Suppose $e:M\rightarrow \R^k$ is positive and satisfies \eqref{ellipinequal}.  Suppose there exists $C_2, p_0\geq 0$, and a smooth function $\nu$ with $ 0\leq \nu\leq 1$   such that 
\begin{equation}\Phi(R,x) := \frac{1}{R^{m-2+p_0}}\int_{B_R(x_0)}(1+\nu(r, x))e\, d\textup{vol}\end{equation} is monotonically increasing in $R$. 
There exists a constant $\epsilon_0>0$, such that if  
\begin{equation}
\Phi(R,x_0) < \epsilon_0,
\end{equation} 
then for any $\delta< \dsst \frac{3}{4}$,
	\begin{equation}
	\sup_{x\in B_{\delta R}(x_0)} e\leq \frac{C}{(\delta R)^{2-p_0}},
	\end{equation}
	where $C$ is an absolute constant only depending on $M$ and the monotonicity constant. 
	
	Moreover, we have
	\begin{equation}\label{fixedsigma}
	\sup_{B_{\frac{1}{2}R}(x_0)} e(f) \leq CR^{-m}\int_{B_{R}(x_0)}(1+\nu(r,x))e(f)\, d\textup{vol}.
	\end{equation}
\end{theorem}

\begin{proof}
	From the monotonicity formula, which applies to all $x\in B_R(x_0)$, we have
	\begin{equation}
	\rho^{2-m-p_0}\int_{B_{\rho}(x)} (1+\nu(r,x))e(f)\,d\textrm{vol} \leq C\tau^{2-m-p_0}\int_{B_{\tau}(x)} (1+\nu(r,x))e(f)\, d\textrm{vol},
	\end{equation}
	for any $x\in B_R(x_0)$ and $0< \rho<\tau<R-|x-x_0|$, where $C_1$ only depends on $m$ and $M$. 
	
	Now we set $\dsst R_1=\delta R=\frac{3}{4}R$. Consider $0<\rho<R_1-|x-x_0|<R-|x-x_0|$. By monotonicity and the fact that $B_{R-|x-x_0|}(x)\subset B_R(x_0)$ the above inequality implies, for any $x\in B_{R_1}(x_0)$, that
	\begin{equation}
	\begin{split}
    &	\rho^{2-m-p_0}\int_{B_\rho(x)}(1+\nu(r,x)) e(u)\, d\textrm{vol}\\ \leq & C_1\left(R_1-|x-x_0|\right)^{2-m-p_0}\int_{B_{R_1-|x-x_0|}(x)}(1+\nu(r,x)) e(f) \, d\textrm{vol}\\
	\leq & C_2\left(R-|x-x_0|\right)^{2-m-p_0} \int_{B_{R-|x-x_0|}(x)}(1+\nu(r,x)) e(f)\, d\textrm{vol} \\
	\leq & C_2\left(\frac{1}{3}R_1\right)^{2-m-p_0} \int_{B_R(x_0)}(1+\nu(r,x)) e(f)\, d\textrm{vol}\\
	= & C\left(\frac{1}{4}\right)^{2-m-p_0}R^{2-m-p_0}\int_{B_R(x_0)}(1+\nu(r,x)) e(f)\, d\textrm{vol}. 
	\end{split}
	\end{equation}

	The last step holds because $R-|x-x_0| = \frac{4}{3}R_1-|x-x_0|\geq \frac{1}{3}R_1$. \\
	
	Now consider the following function 
	\begin{equation}
	h(\sigma)= (R_1-\sigma)^{2-p_0}\sup_{B_\sigma(x_0)} e(f).
	\end{equation}
	
	Let $\sigma_0$ be the point where $h(\sigma_0)=\dsst \max_{0\leq \sigma < R_1}h(\sigma)$. Suppose at $x_1\in \overline{B_{\sigma_0}(x_0)}$ we have 
	\begin{equation}
	e(f)(x_1) = B_{\sigma_0}(x_0) e(f) := e_0.
	\end{equation}
	
	Let $\rho_0=\dsst \frac{1}{2}(R_1-\sigma_0)$. Since $B_{\rho_0}(x_1)\subset B_{\rho_0+\sigma_0}(x_0)$, we have
	\begin{equation}
	\sup_{B_{\rho_0}(x_1)} e(f)\leq \sup_{B_{\rho_0+\sigma_0}(x_0)} e(f).
	\end{equation}

	From the definition of $h(\sigma)$, we get
	\begin{equation}
	(B_1-(\rho_0+\sigma_0))^2\sup_{B_{\rho_0+\sigma_0}(x_0)} e(f) \leq (B_1-\sigma_0)^2\sup_{B_{\sigma_0}(x_0)} e(f). 
	\end{equation}
	
	This implies
	\begin{equation}
	\sup_{B_{\rho_0}(x_1)} e(f)\leq \sup_{B_{\rho_0+\sigma_0}(x_0)} e(f) \leq 4\sup_{B_{\sigma_0}(x_0)} e(f) = 4e_0.
	\end{equation}
	
	Now consider the elliptic equation in $B_{\rho_0}(x_1)$. We have
	\begin{equation}
	\Delta e(f) \leq Ce^2(f) \leq C\left(\sup_{B_{\rho_0}(x_1)} e(f)\right)e(f)\leq 4C e_0 e(f)\qquad \text{in }B_{\rho_0}(x_1).
	\end{equation}
	
	By the single point Moser iteration result, we have
	\begin{equation}
	e_0 = e(x_1) \leq \left(C_1(M)4Ce_0+\frac{C_2(M)}{\rho_0^2}\right)^{\frac{m}{2}}\int_{B_{\rho_0}(x_1)} e(f)\, d\textrm{vol}.
	\end{equation}

	Rearranging, we get
	\begin{equation}
	\begin{split}
	\rho_0^{2-p_0}e_0 \leq & \left(C_1(M)\rho_0^2e_0+C_2(M)\right)^{\frac{m}{2}}\frac{1}{\rho_0^{m-2+p_0}}\int_{B_{\rho_0}(x_1)} e(f)\, d\textrm{vol} \\
	\leq & \left(C_1(M)\rho_0^{2-p_0}e_0+C_2(M)\right)^{\frac{m}{2}}\frac{1}{\rho_0^{m-2+p_0}}\int_{B_{\rho_0}(x_1)} (1+\nu(r,x))e(f)\, d\textrm{vol}.
	\end{split}
	\end{equation}
	
	This is
	\begin{equation}
	h(\sigma_0)\leq \left(C_1(M)h(\sigma_0)+C_2(M)\right)^{\frac{m}{2}}\frac{1}{\rho_0^{m-2+p_0}}\int_{B_{\rho_0}(x_1)}(1+\nu(r,x)) e(f)\, d\textrm{vol}.
	\end{equation}
	
	Since $x_1\in B_{R_1}(x_0)$ and $\rho_0+\sigma_0<R_1=\frac{3}{4}R$, we can apply the previous result that
	\begin{equation}
	\frac{1}{\rho_0^{m-2+p_0}}\int_{B_{\rho_0}(x_1)}(1+\nu(r,x)) e(f)\, d\textrm{vol}\leq C\frac{1}{R^{m-2+p_0}}\int_{B_R(x_0)}(1+\nu(r,x)) e(f)\, d\textrm{vol},
	\end{equation}
	and therefore,
	\begin{equation}\label{inequaForH}
	h(\sigma_0) \leq \left(C_1(M)h(\sigma_0)+C_2(M)\right)^{\frac{m}{2}}\frac{1}{R^{m-2+p_0}}\int_{B_R(x_0)}(1+\nu(r,x)) e(f)\, d\textrm{vol}.
	\end{equation}
	
	Suppose there exists $R_1=\frac{3}{4}R$ such that $h(\sigma_0)=1$. Then
	\begin{equation}
	1 = h(\sigma_0) \leq \left(C_1+C_2\right)^{m/2}\frac{1}{R^{m-2+p_0}}\int_{B_R(x_0)}(1+\nu(r,x)) e(f)\, d\textrm{vol}.
	\end{equation}
	
	We can pick $\epsilon_0\leq \dsst \frac{1}{2(C_1+C_2)^{m/2}}$ such that 
	\begin{equation}
	\left(C_1+C_2 \right)^{m/2}\frac{1}{R^{m-2+p_0}}\int_{B_R(x_0)} (1+\nu(r,x))e(f)\, d\textrm{vol} \leq (C_1+C_2)^{m/2}\epsilon_0 <1, 
	\end{equation}
	for any $R$ such that $\Psi(R, x_0)\leq \epsilon_0$. This  leads to a contradiction. Therefore, \\$h(\sigma_0)\leq 1$. \\
	
	That is,
	\begin{equation}
	h(\sigma) = (R_1-\sigma)^{2-p_0}\sup_{B_\sigma(x_0)}e(f) \leq h(\sigma_0) \leq 1.
	\end{equation}
	
	We can see that the above proof will hold when we replace $R_1=\frac{3}{4}R$ with $R_1=\delta R$ with any $\delta\leq \dsst \frac{3}{4}$. \\
	
	Finally, picking $\sigma=\frac{1}{2}R_1=\frac{1}{2}\delta R$, we get
	\begin{equation}
	\sup_{B_{\frac{1}{2}\delta R}(x_0)} e(f) \leq \frac{C}{(\frac{1}{2}\delta R)^{2-p_0}}.
	\end{equation}
	with $C$ a constant depends on $C_0$. (We can denote new $\delta$ as $\frac{1}{2}\delta$ to get the desired result.)\\
	
	As for \eqref{fixedsigma}, we go back to the original estimate \eqref{inequaForH} and use the fact that $h(\sigma_0)\leq 1,$
	\begin{equation}
	\begin{split}
	h(\sigma_0)\leq &  \left(C_1(M)\rho_0^{p_0}+C_2(M)\right)^{\frac{m}{2}}\frac{1}{\rho_0^{m-2+p_0}}\int_{B_{\rho_0}(x_1)}(1+\nu(r,x)) e(f)\, d\textrm{vol}\\
	\leq & C\frac{1}{R^{m-2+p_0}}\int_{B_{R}(x_0)} (1+\nu(r,x))e(f)\, d\textrm{vol}.
	\end{split}
	\end{equation}
	
	Using the definition of $h(\sigma_0)$ again, we have
	\begin{equation}
	(R_1-\sigma)^{2-p_0}\sup_{B_{\sigma}(x_0)}(1+\nu(r,x)) e(f) \leq h(\sigma_0) \leq C\frac{1}{R^{m-2+p_0}}\int_{B_{R}(x_0)} e(f)\, d\textrm{vol}.
	\end{equation}
	
	Here we pick $\delta$ to be $\frac{3}{4}$, and $\sigma = \frac{1}{2}R<R_1=\delta R$, we have
	\begin{equation}
	\sup_{B_{\frac{1}{2}R}(x_0)} e(f) \leq C R^{-m}\int_{B_{R}(x_0)} (1+\nu(R,x))e(f)\, d\textrm{vol}.
	\end{equation}
	
\end{proof}




\section{Parabolic $\epsilon$-Regularity}

In this section, we consider a positive function $e:M\rightarrow\R^1$ that satisfies the following parabolic inequality
\begin{equation}\label{parainequal}
\frac{\partial e}{\partial t}+\Delta e\leq Ce^2 \qquad \text{in }P_R(z_0)
\end{equation}
for some $R<R_M$. Here $z_0=(x_0, t_0)$. The regions $P_R(z_0)$ and $T_R(t_0)$ are the same as the ones defined in Chapter 2.


\begin{theorem}
	Suppose $e:M\times I\rightarrow \R^k$ is a positive function that satisfies \eqref{parainequal}. Suppose there exists $ \nu_0\geq 0$ such that 
	\begin{equation}
	\Psi((x_0,t_0), R) := \int_{t_0-4R^2}^{t_0-R^2}\frac{1}{(t_0-t)^{\nu_0}}\int_M G_{(x_0,t_0)}(x,t)e(f)(x,t)\varphi^2\, d\textup{vol}\, dt  \end{equation} 
	is monotonically increasing in $R$. There exists $\epsilon_0$ such that, if there exists $R\leq \sqrt{t_0-t}/2$ with
	\begin{equation}
	\Psi((x_0, t_0), R)=\int_{t_0-4R^2}^{t_0-R^2}\frac{1}{(t_0-t)^{\nu_0}}\int_{M}G_{(x_0,t_0)}(x,t)e(x,t)\varphi^2\, d\textup{vol}\, dt \leq \epsilon_0,
	\end{equation}
     then for any $\delta$ that depends only on $M, E(f_0(x)), R$, and $C$ an absolute constant, we have
     \begin{equation}
     \sup_{P_{\delta R}(x_0, t_0)}e\leq \frac{C}{R^{2-2\nu_0}}.
     \end{equation}
	
	Here $G_{(x_0, t_0)}$ is the backward Gaussian 
	\begin{equation}
	G_{(x_0, t_0)}(x,t)=\frac{1}{\left(4\pi(t_0-t)\right)^{\frac{m}{2}}}\exp\left(-\frac{|x-x_0|^2}{4(t_0-t)}\right),
	\end{equation}
	and $\varphi(x)=\varphi(|x-x_0|)$ is a cut-off function supported in $B_{R_M}(x_0)$. (Sometimes we denote $G_{(x_0, t_0)}(x,t)=G(x,t)$ or just $G_{(x_0,t_0)}$, $G_{z_0}$.)	
\end{theorem}

\begin{proof}
	Let $R_1=\delta R$ for $\delta <1/2$. For any $\sigma <R_1$, consider the following situation. Suppose
	$$ e_0 = \sup_{P_\sigma (z_0)} e(f). $$
	
	We are interested in the point $z_1=(x_1,t_1) \in \overline{P_\sigma(z_0)}$ such that 
	$$ e(f)(z_1) = e_0. $$
	
	Then inside some cylinder $P_{\rho}(z_1)\subset P_{R_1}(z_0)$, we have
	$$ \left(\frac{\partial}{\partial t}+\Delta\right)e(f) \leq C\left(\sup_{P_{R_1-\sigma}(z_0)} e(f)\right) e. $$

	From the single point Moser iteration, we have
	$$ e(f)(z_1) = \sup_{P_{\sigma}(z_0)} e(f) \leq \left(C_1(M)\left(\sup_{P_{\rho}(z_1)} e(f)\right)+\frac{C_2(M)}{\rho^2}\right)^{\frac{m+2}{2}}\int_{P_\rho(z_1)} e(f)d\textrm{vol}. $$
	
	Rearranging the above inequality so that it attains a form in which we may apply the monotonicity formula, we get
	\begin{equation*}
	\begin{split}
	& \rho^{2-2\nu_0}\sup_{P_{\sigma}(z_0)} e(f)\\
	\leq & \left(C_1(M)\rho^{2-p_0}(\sup_{P_{\rho}(z_1)}e)\rho^{p_0} + C_2(M)\right)^{(m+2)/2}\frac{1}{\rho^{m+2\nu_0}} \int_{P_\rho(z_1)} e(f)\, d\textrm{vol}\, dt.
	\end{split}
	\end{equation*} 
	
	To relate $\dsst \sup_{P_\rho(z_1)}e(f)$ back to $\dsst \sup_{P_{\sigma}(z_0)}e(f)$, we define the function
	\begin{equation}
	h(\sigma) = (R_1-\sigma)^{2-2\nu_0}\sup_{P_\sigma(z_0)} e(f),
	\end{equation}
	and pick $\rho = \frac{1}{2}(R_1-\sigma)$. 
	
	Since $e(f)$ is regular in $P_R(z_0)$, we are able to find a point $\sigma_0$ such that 
	$$ h(\sigma_0) = \max_{0\leq \sigma < R_1}h(\sigma). $$
	
	There exists $z_1=(x_1,t_1) \in \overline{P_{\sigma_0}(z_0)}$ such that
	$$ e(f)(z_1) = \sup_{P_{\sigma_0}(z_0)} e(f). $$
	
	Let $\dsst \rho_0 = \frac{1}{2}(R_1-\sigma_0)$. We will argue that 
	$$ \sup_{P_{\rho_0}(z_1)} e(f) \leq 4 \sup_{P_{\sigma_0}(z_0)} e(f). $$
	
	It is easy to check that $\rho_0+\sigma_0=\dsst \frac{1}{2}(R_1+\sigma_0)<R_1=\delta R<R$ and $P_{\rho_0}(z_1)\subset P_{\rho_0+\sigma_0}(z_0)$. Similar argument as in the elliptic case, we have
	\begin{equation}
	\sup_{P_{\rho_0}(z_1)} e(f) \leq \sup_{P_{\rho_0+\sigma_0}(z_0)} \leq 4\sup_{P_{\sigma_0}(z_0)} = 4e(f)(z_1).
	\end{equation}
	
	Now we base our estimate inside the cylinder $P_{\rho_0}(z_1)$ since $z_1$ is the point that $e(f)$ achieves its max within $P_{\sigma_0}(z_0)$. From the single point Moser iteration, we get
	\begin{equation}
	\begin{split}
	e(f)(z_1) \leq &  \left(C_1(M)\left(\sup_{P_{\rho_0}(z_1)}e(f)\right)+\frac{C_2(M)}{\rho_0^2}\right)^{\frac{m+2}{2}}\int_{P_{\rho_0}(z_1)} e(f)\, d\textrm{vol}\, dt\\
	\leq & \left(4C_1(M)\left(\sup_{P_{\sigma_0}(z_0)}e(f)\right)+\frac{C_2(M)}{\rho_0^2}\right)^{\frac{m+2}{2}}\int_{P_{\rho_0}(z_1)} e(f)\, d\textrm{vol}\, dt.
	\end{split}
	\end{equation}
	
	Do the same rearranging as before, we get the following,
	\begin{equation}
	\begin{split}
	& \rho_0^{2-2\nu_0} e(f)(z_1)\\
	 = & \frac{1}{2^{2-2\nu_0}}(R_1-\sigma_0)^{2-2\nu_0}\sup_{P_{\sigma_0}(z_0)} e(f)\\
	\leq & \left(C_1(M)(R_1-\sigma_0)^{2-2\nu_0}\sup_{P_{\sigma_0}(z_0)} e(f)+C_2(M)\right)^{(m+2)/2}\frac{1}{\rho_0^{m+2\nu_0}} \int_{P_{\rho_0}(z_1)} e(f)\,d\textrm{vol}\,dt.
	\end{split}
	\end{equation}
	
	That is,
	$$ h(\sigma_0) \leq \left(C_1(M)h(\sigma_0)+C_2(M)\right)^{(m+2)/2}\frac{1}{\rho_0^{m+2\nu_0}} \int_{P_{\rho_0}(z_1)} e(f)\, d\textrm{vol}\,dt. $$
	
	Here, as long as we can show that 
	$$ \frac{1}{\rho_0^{m+2\nu_0}} \int_{P_{\rho_0}(z_1)} e(f)\, d\textrm{vol}\,dt $$
	can be bounded by $\epsilon_0$ using the monotonicity formula, we could argue that $h(\sigma_0)$ cannot be arbitrarily large, provided $\sigma_0<R_1=\delta R<R$, where $\epsilon_0$ gets achieved. Therefore, we will have control for all $\dsst \sup_{P_{\sigma}(z_0)} e(f),\quad \forall \sigma<R_1$ from $h(\sigma_0)$. 
	
	We present the following lemma and will prove it after the proof of this theorem.
	\begin{lemma}\label{paraepsilonlemma}
		Let $R_1=\delta R$ for some $\delta<\frac{1}{2}$, and $R<R_M$. Let $\rho,\ \sigma \in(0, R_1)$ with $\rho+\sigma <R_1$. Then for any $z_1=(x_1,t_1)\in P_\sigma(z_0)$, the monotonicity formula implies
		\begin{equation}
		\begin{split}
		& \frac{1}{\rho^{m+2\nu_0}} \int_{P_\rho(z_1)} e(u)\, d\textup{vol}\, dt\\
		\leq & C_1(R)\Psi((x_0,t_0), R)+C_2RE(f_0(x))+C_3\exp\left((2-m-2\nu_0)\ln R-c_2\delta^{-2}\right)E(f_0(x)).
		\end{split}
		\end{equation}
		Here $C_1(R)$ is a constant depends on $R$ and $M$, $C_3$, $c_2$ are absolute constants. $E(f_0(x))$ is the initial energy.
	\end{lemma}
	
	With the lemma, we can proceed to argue for bounds of $h(\sigma_0)$ as follows. Assume $\Psi((x_0,t_0),R_0)\leq \epsilon_0$. We then pick $\delta$ small so that\\ $C_3\exp\left((2-m-2\nu_0)\ln  R-c_2\delta^{-2}\right)E(f_0(x)) \leq \epsilon_0$, for some small $R<R_M$. We have
	\begin{equation}
	\begin{split}
	h(\sigma_0) \leq & \left(C_1h(\sigma_0)+C_2\right)^{(m+2)/2}\frac{1}{\rho_0^{m+2\nu_0}}\int_{P_{\rho_0}(z_1)} e(f)\, d\textrm{vol}\, dt\\
	\leq & \left(C_1h(\sigma_0)+C_2\right)^{(m+2)/2}(c_1\epsilon_0+c_2RE(f_0(x))) \\
	\leq & \left(C_1h(\sigma_0)+C_2\right)^{(m+2)/2}\epsilon_0.
\end{split}	
\end{equation}
	
	Suppose there exists $R_1=\delta R$, such that for any $\sigma<R_1$, $h(\sigma_0)=\dsst \max_{0\leq \sigma < R_1}=1$. By the above inequality, we have
	\begin{equation}
	1 \leq \left(C_1+C_2\right)^{(m+2)/2}\epsilon_0.
	\end{equation} 	
	
	Therefore, we can choose $\epsilon_0$ so that $\left(C_1+C_2\right)^{(m+2)/2}\epsilon_0<1$. That said, with $R_1=\delta R<\delta R_0$, and for all $\sigma<R_1$, we always have
	\begin{equation}
	h(\sigma)\leq h(\sigma_0) \leq 1.
	\end{equation} 
	
	Since $h(\sigma_0)=\dsst \max_{0\leq \sigma < R_1} h(\sigma)$, we can pick $\sigma = \dsst \frac{1}{2}R_1=\dsst \frac{1}{2}\delta R$, and we have
	\begin{equation}
	(R_1-\sigma)^{2-p_0}\sup_{P_\sigma(z_0)} e(f) \leq h(\sigma_0)  \leq 1,
	\end{equation} 
	which is,
	\begin{equation}
	\sup_{P_{\frac{1}{2}\delta R}(z_0)} e(f) \leq \frac{C}{(\frac{1}{2}\delta R)^{2-p_0}},
	\end{equation}
	where $C$ is an absolute constant. Picking new $\delta$ to be half of the original $\delta$, we will get the desired estimate.

	
	
	
	
\end{proof}

\begin{proof}[Proof of Lemma \ref{paraepsilonlemma}]
	
	We want to show that
	$$ \frac{1}{\rho^{m+2\nu_0}}\int_{P_\rho(z_1)}e(f)\, d\textrm{vol}\, dt $$
	can be controlled by some $\Psi((x_0,t_0), R)$ and eventually bounded by $\epsilon_0$. \\
	
	First, we need an estimate between the backward Gaussian $G$ and $\rho^{-m}$. Since Gaussians blow up at $t=t_0$, we cannot derive the inequality directly from $z_1 = (x_1, t_1)$. To get a bound on the Gaussian, we move the center a little above so that it does not lie in the region $T_R(t_0)=M\times [t_0-4R^2, t_0-R^2]$, where we later apply the monotonicity formula. \\
	
	Let $R_1 = \delta R$ and $\sigma, \rho \in(0, R_1)$ such that $\rho+\sigma <R_1$. Let $z_1=(x_1, t_0)\in P_\sigma(z_0)$. For any point in $z\in P_\rho(z_1)$, we have the following estimate 
	\begin{equation}
	\begin{split}
	G_{(x_1, t_1+2\rho^2)} = & (4\pi(t_1+2\rho^2-t))^{-m/2}\exp\left(-\frac{|x-x_1|^2}{4(t_1+2\rho^2-t)}\right)\\
	\geq & \left(4\pi(t_1+2\rho^2-(t_1-\rho^2))\right)^{-m/2}\exp\left(-\frac{\rho^2}{4(t_1+2\rho^2-(t_1+\rho^2))}\right)\\
	= & C \rho^{-m}.
	\end{split}
	\end{equation}
	
	Therefore, we have
	\begin{equation}
	\rho^{-m-2\nu_0}\int_{P_\rho(z_1)} e(f)\, d\textrm{vol}\, dt \leq C_1\rho^{-2\nu_0}\int_{P_\rho(z_1)} G_{(x_1, t_1+2\rho^2)}(x,t)e(f)\, d\textrm{vol}\, dt.
	\end{equation}
	
	It is obvious that $P_\rho(z_1)\subset T_\rho(t_1+2\rho^2)$. Suppose $\varphi^2$ is the cut-off function that is used in the proof of the monotonicity formula that is compactly supported in $B_{R_M}(x_0)$, we have
	\begin{equation}
	\begin{split}
	& \rho^{-2\nu_0}\int_{P_\rho(z_1)} G_{(x_1, t_1+2\rho^2)}(x,t)e(f)\, d\textrm{vol}\, dt\\
	\leq &\rho^{-2\nu_0}\int_{T_\rho(t_1+2\rho^2)} G_{(x_1, t_1+2\rho^2)}(x,t)e(f)\varphi^2\, d\textrm{vol}\, dt \\
	\leq & C\int_{t_1-3\rho^2}^{t_1+\rho^2}\frac{1}{(t_1+2\rho^2-t)^{\nu_0}}\int_M G_{(x_1, t_1+2\rho^2)}(x,t)e(f)\varphi^2\, d\textrm{vol}\, dt\\
	= & C\Psi((x_1, t_1+2\rho^2), \rho).
	\end{split}
	\end{equation}
	
	By the monotonicity formula, we have
	\begin{equation}
	\Psi((x_1, t_1+2\rho^2), \rho)\leq C_1(R)\Psi((x_1, t_1+2\rho^2), \sqrt{3}R/2)+C_2RE(f_0(x)).
	\end{equation}
	
	With the fact that $R_1=\delta R\leq \frac{1}{2}R$ and $\rho,\sigma<R_1$, we can check that $t_1+2\rho^2-3R^2>t_0-4R^2$. We then split the stripe $T_{R/2}(t_1+2\rho^2)$ into 3 parts
	\begin{equation}
	\begin{split}
	& \Psi((x_1, t_1+2\rho^2), \sqrt{3}R/2)\\
	\leq & \left(\int_{t_0-4R^2}^{t_0-R^2} +\int_{t_0-R^2}^{t_1+2\rho^2-3R^2/4}\right)\frac{1}{(t_1+2\rho^2-t)^{\nu_0}}\int_M G_{(x_1, t_1+2\rho^2)}(x,t)e(f)\varphi^2\, d\textrm{vol}\, dt\\
	\leq & C\int_{T_R(t_0)}\frac{1}{(t_0-t)^{\nu_0}}G_{(x_1,t_1+2\rho^2)}e(f)\varphi^2\, d\textrm{vol}\, dt\\
	&  + \int_{t_0-R^2}^{t_1+2\rho^2-3R^2/4}\frac{1}{(t_1+2\rho^2-t)^{\nu_0}}\int_M G_{(x_1,t_1+2\rho^2)}e(f)\varphi^2\, d\textrm{vol}\, dt\\
	= & \text{I} + \text{II}.
	\end{split}
	\end{equation}
	
	For $\text{II}$, since $\sigma^2+\rho^2\leq (\sigma+\rho)^2<R_1^2<R^2/4$, we apply the monotonicity formula again and get
	\begin{equation}
	\begin{split}
	& \int_{t_0-R^2}^{t_1+2\rho^2-3R^2/4}\frac{1}{(t_1+2\rho^2-t)^{\nu_0}}\int_M G_{(x_1, t_1+2\rho^2)}e(f)\varphi^2\,d\textrm{vol}\, dt \\
	\leq &   C\int_{t_0-R^2}^{t_0+\sigma^2+2\rho^2-3R^2/4}\frac{1}{(t_0-t)^{\nu_0}}G_{(x_1, t_1+2\rho^2)}e(f)\varphi^2\, d\textrm{vol}\,dt\\
	\leq & C\int_{t_0-R^2}^{t_0-R^2/4}\frac{1}{(t_0-t)^{\nu_0}}G_{(x_1, t_1+2\rho^2)}e(f)\varphi^2\, d\textrm{vol}\,dt\\
	\leq & C\int_{t_0-4R^2}^{t_0-R^2}\frac{1}{(t_0-t)^{\nu_0}}G_{(x_1, t_1+2\rho^2)}e(f)\varphi^2\, d\textrm{vol}\,dt.
	\end{split}
	\end{equation}
	
	So far we have shown
	\begin{equation}
	\begin{split}
	 & \rho^{-m-2\nu_0}\int_{P_\rho(z_1)} e(f)\, d\textrm{vol}\, dt\\
	  \leq & C_1(R)\int_{T_R(t_0)} \frac{1}{(t_0-t)^{\nu_0}}\int_M G_{(x_1, t_1+2\rho^2)}(x,t)e(f)\varphi^2\, d\textrm{vol}\, dt + C_2(R)E(f_0(x)).
	  \end{split}
	\end{equation}
	
	Now we need to get a bound for $G_{(x_1, t_1+2\rho^2)}(x,t)$ from $G_{x_0, t_0}$. Here we consider two regions: $|x-x_1|\geq R/\delta$ and $|x-x_1|<R/\delta$.\\ 
	
	\textbf{Case 1}: when $|x-x_1|\geq R/\delta$.\\
	
	For $(x,t)\in T_R(t_0)$, $t_0-t\geq R$ and $t_1-t_0+2\rho^2\geq -\sigma^2\geq -\frac{1}{2}R^2\geq -\frac{1}{2}(t_0-t)$. Therefore, $(t_1+2\rho^2-t)^{m/2}\geq c(t_0-t)^{m/2}$. Hence, 
	\begin{equation}
	\begin{split}
	G_{(x_1, t_1+2\rho^2)}(x,t) = & \left(4\pi(t_1+2\rho^2-t)\right)^{-\frac{m}{2}}\exp\left(-\frac{|x-x_1|^2}{4(t_1+2\rho^2-t)}\right)\\
	\leq & c\left(4\pi(t_0-t)\right)^{-\frac{m}{2}}\exp\left(-\frac{(R/\delta)^2}{4R^2+\sigma^2+2\rho^2}\right) \\
	\leq & c_1R^{-m}\exp\left(-c_2\delta^{-2}\right).
	\end{split}
	\end{equation}
	
	\textbf{Case 2}: when $|x-x_1|<R/\delta$. \\
	
	We use the same argument that \\$t_1+2\rho^2-t\leq t_0-t+\sigma^2+2\rho^2\leq t_0-t+2\delta^2R^2\leq \left(1+2\delta^2\right)(t_0-t)$, and the triangle inequality to get
	\begin{equation}
	\begin{split}
	& G_{(x_1, t_1+2\rho^2)}(x,t) \\
	\leq & c\left(4\pi(t_0-t)\right)^{-\frac{m}{2}}\exp\left(-\frac{|x-x_0|^2}{4(t_0-t)}\right)\times \exp\left(\frac{|x-x_0|^2}{4(t_0-t)}-\frac{|x-x_1|^2}{4(t_1+2\rho^2-t)}\right)\\
	= &  cG_{(x_0, t_0)}\exp\left(\frac{|x-x_0|^2}{4(t_0-t)}-\frac{|x-x_1|^2}{4(t_1+2\rho^2-t)}\right)\\
	\leq & cG_{(x_0,t_0)}\exp\left(\frac{1}{4(t_0-t)}\left[|x-x_0|^2-\frac{|x-x_1|^2}{1+2\delta^2}\right]\right)\\
	\leq & cG_{(x_0,t_0)}\exp\left(\frac{1}{4(t_0-t)}\left[\left(|x-x_1|+|x_1-x_0|\right)^2-\frac{|x-x_1|^2}{1+2\delta^2}\right]\right)\\
	= & cG_{(x_0,t_0)}\exp\left(\frac{1}{4(t_0-t)}\left[\frac{2\delta^2}{1+2\delta^2}|x-x_1|^2+2|x-x_1||x-x_0|+|x-x_0|^2\right]\right)\\
	\leq & cG_{(x_0, t_0)}\exp\left(\frac{1}{4(t_0-t)}2\delta^2\left(\frac{R}{\delta}\right)^2+2\left(\frac{R}{\delta}\right)\sigma+\sigma^2\right)\\
	\leq & cG_{(x_0, t_0)}\exp\left(\frac{1}{4(t_0-t)}c_2R^2\right)\\
	\leq & cG_{(x_0, t_0)}\exp\left(c_2/4\right) = \tilde{c}G_{(x_0, t_0)}.
	\end{split}
	\end{equation}
	
	Therefore, we have for any $(x,t)\in T_R(t_0)$,
	\begin{equation}
	\begin{split}
	G_{(x_1,t_1+2\rho^2)}(x,t)\leq & cG_{(x_0,t_0)}+c_1R^{-m}\exp\left(-c_2\delta^{-2}\right) \\
	\leq & cG_{(x_0, t_0)}(x,t)+c_1R^{-2+2\nu_0}\exp\left((2-m-2\nu_0)\ln R-c_2\delta^{-2}\right).
	\end{split}
	\end{equation}
	
	With the estimate of the Gaussian, we have for any $z_1\in P_\sigma(z_0)$, and \\$\sigma+\rho<R_1=\delta R<R/2$, 
	\begin{equation}
	\begin{split}
	& \rho^{-m-2\nu_0}\int_{P_\rho(z_1)} e(f)\, d\textrm{vol}dt \\
	\leq & C_1(R)\int_{T_R(t_0)}\frac{1}{(t_0-t)^{2\nu_0}}\int_M \left[\left(cG_{(x_0, t_0)}+c_1R^{-2+2\nu_0}\exp\left((2-m-2\nu_0)\ln R-c_2\delta^{-2}\right)\right)\right.\\
	& \times \left. e(f)\varphi^2\right]\, d\textrm{vol}\, dt + C_2(R)E(f_0(x)) \\
	\leq & C_1(R)\Psi((x_0,t_0), R)+\left[C_3(R)\exp\left((2-m-2\nu_0)\ln R-c_2\delta^{-2}\right)+C_2(R)\right]E(f_0(x))\\
	= & C_1(R)\Psi((x_0, t_0), R)+C_3(R)\exp\left((2-m-2\nu_0)\ln R-c_2\delta^{-2}\right)E(f_0(x))+C_2RE(f_0(x)).
	\end{split}
	\end{equation}
	
	To get the second term small, we may pick $\delta \approx c(m,\nu_0)|\ln R|^{-1/2}$ so that $c_2\delta^{-2}$ is much larger than $(2-m-2\nu_0)\ln R$ when $R<1$, and $\delta$ can be chosen independently if $R\geq 1$. 

\end{proof}

%% file: Chapter5/Chapter5.tex
\chapter{Models with Repulsive Potentials}

\section{Sigma Models with Repulsive Potentials}

The sigma model approach from \cite{StruweSigma} illustrates an idea to replace the nonlinear constraints from the target manifold with potential energy constraints. This approach allows us to replace a nonlinear problem that may not have solutions with a sequence of problems each of which has a solution. Then the problem becomes showing this sequence of solutions converges to a solution of the original problem.

The repulsive potential is motivated by the desire to keep the homotopy class of the projection of $f$ onto the Grassmannian fixed. We are particularly interested in the existence of harmonic maps from a compact Riemannian manifold $M$ to complex Grassmannian manifolds. It is well-known that Grassmannian manifolds have positive sectional curvature. So the approach from \cite{EellsSampson} cannot be applied. However, Grassmannian manifolds with their canonical metrics are Kahler, and holomorphic maps between Kahler manifolds are harmonic. Hence it is easy to show there exist many harmonic maps to the Grassmannian. The difficulty is producing them without knowing holomorphicity in advance.  

The rank $k$ complex Grassmannian $\text{Gr}(k, \mathds{C}^l)$ can be defined as the set of Hermitian matrices in $\R^{l\times l}$ that have $k$ $(+1)$ eigenvalues and the remaining eigenvalues are 0. An equivalent representation is:
\begin{equation}
\text{Gr}(k, \mathds{C}^l) = \left\{u \in \mathds{R}^{l\times l}|\text{$u$ is hermitian and $k$ eigenvalues of $u$ are 1, the rest are -1} \right\}.
\end{equation}  

Consider the classifying map $f_F:M\rightarrow \text{Gr}(k, \mathds{C}^l)$ of a vector bundle $F$ over $M$, and consider the energy
$$ E(f_F) = \int_M \left(\frac{1}{2}|df_F|^2+W(f_F)\right)\, d\textrm{vol}. $$

To guarantee that the projection of energy minimizers to the
Grassmannian remain in the same homotopy class, we have introduced a potential term $W$ to our energy. $W(f_F)$ is chosen so that $W(f_F)\rightarrow \infty$ if one of the eigenvalues of $f_F\rightarrow 0$, and $W$ is smooth in the complement of the set $\left\{ \det(f_F)=0\right\}$.

Since $W(f_F)$ is smooth as long as all eigenvalues of $f_F$ are away from 0, singularities of $f_F$ can only arise when some eigenvalues of $f_F$ approach to 0. This translates questions about singularities of $f_F$ into questions about singularities of $W(f_F)$.

Because $W(f_F)$ is pushing $f_F$ away from the set $\left\{\det(f_F)=0\right\}$, we call $W$ a ``repulsive potential".

\section{Repulsive Potentials}

Let $f:M\rightarrow R^{l\times l}\supset \text{Gr}(k, \mathds{C}^l)$. Let $\{\lambda_i\}_{i=1}^l$ be the eigenvalues of $f$.  Consider the energy functional
$$ E(f) = \int_M \left(\frac{1}{2}|df|^2+W(f)\right)\, d\textrm{vol}, $$
with repulsive potential 
\begin{equation} 
W: R^{l\times l}\rightarrow R^1, \qquad \text{with}\qquad W(f)=\frac{1}{2}\text{trace}(f^{-2})=\frac{1}{2}\sum_i^l \frac{1}{\lambda_i^2}:= \frac{1}{2}|f^{-1}|^2.
\end{equation}
Here $\left|\,\cdot\, \right|$ is the Hilbert-Schmidt norm of matrices.\\

When the eigenvalues of $f$ are bounded, we have bounded potential $W(f)$ and $f$ lies in a set that is homotopy equivalent to the Grassmannian $\textrm{Gr}(k, \mathds{C}^l)$.

This energy functional $E$ has a nice scaling property, if we replace $M$ with a Euclidean space. Assuming
$$ f = \lambda \tilde{f},\qquad x = \lambda^2 y, $$
we have
\begin{equation}
\begin{split}
E(\tilde{f}(y)) = &\int_M\left( \frac{1}{2}|d_y\tilde{f}|^2+\frac{1}{2}|\tilde{f}^{-2}|\right)\, d\textrm{vol}_y\\
 = &\frac{1}{\lambda^{2m-2}} \int_M \left(\frac{1}{2}|d_xf|^2+\frac{1}{2}|f^{-2}|\right)\, d\textrm{vol}_x\\
  = &\frac{1}{\lambda^{2m-2}}E(f(x)). 
\end{split}
\end{equation} 

\subsection{Estimates of $\nabla W$ and $\text{Hess}(W)$}
For this specific potential $W$, we calculate some of its related quantities, which will be used later. 

The directional derivative of $W$ along the $B$ direction is
\begin{equation}
\begin{split}
W_B(f) = & \frac{d}{ds}W(f+sB)\Big|_{s=0}\\
= & \frac{d}{ds}\frac{1}{2}\text{trace}\left((f+sB)^{-2}\right)\Big|_{s=0}\\
= & \frac{1}{2}\text{trace}\left(\frac{d}{ds}(f+sB)^{-2}\right)\Big|_{s=0}\\
= & \frac{1}{2}\text{trace}\left(-f^{-2}Bf^{-1}-f^{-1}Bf^{-2}\right) = -\text{trace}\left(f^{-3}B\right).
\end{split}
\end{equation}

Therefore
\begin{equation}
\nabla W(f) = -f^{-3},\qquad \text{and}\qquad |\nabla W|^2=|f^{-3}|^2 .
\end{equation}

Under this $W$, the static elliptic equation and the corresponding parabolic gradient equation are
\begin{equation}
d^*df-f^{-3} = 0,
\end{equation}
and
\begin{equation}\label{paraeq}
\frac{\partial f_t}{\partial t} + d^*df_t - f_t^{-3} = 0,\qquad f(0, x) = f_0(x) \in \text{Gr}(k, \mathds{C}^l).
\end{equation}

\subsection{Short Time Existence of Parabolic Flow}

We assume $f_0$, the initial data of \eqref{paraeq}, is smooth. We can convert \eqref{paraeq} into the integral equation
\begin{equation}
f(t, x) = e^{-td^*d}f_0 + \int_0^t e^{(s-t)d^*d}f^{-3}_s\, ds.
\end{equation} 

Consider the subset
\begin{equation}
Z= \left\{f\in C^1(M)|\left\|f-f_0\right\|_{C^1(M)}\leq \alpha \right\},
\end{equation}
for a fixed $\alpha<1$. We have (from \cite[p.272, (1.3) and p.273, (1.11)]{taylor1997partial})
\begin{equation}
e^{-td^*d} : C^1(M)\rightarrow C^1(M) \textrm{ is a strongly continous semigroup for $t\geq 0$},
\end{equation}
and
\begin{equation}
\left\|e^{-td^*d}\right\|_{\mathcal{L}(C^0(M), C^1(M))} \leq Ct^{-1/2},\qquad \textrm{for $t\in (0,1)$}. 
\end{equation}

Furthermore, the functional 
\begin{equation}
F: f\rightarrow f^{-3}
\end{equation}
is defined from $C^1(M)$ to $C^0(M)$ when $f$ is invertible. 

The derivative of $F$ is
\begin{equation}
DF_f(g) = -f^{-1}gf^{-3}-f^{-2}gf^{-2}-f^{-3}gf^{-1}.
\end{equation}

We have the estimate
\begin{equation}
\begin{split}
\left\|DF_f(g)\right\|_{C^0(M)} \leq & \left\|f^{-1}gf^{-3}\right\|_{C^0} + \left\|f^{-2}gf^{-2}\right\|_{C^0}+\left\|f^{-3}gf^{-1}\right\|_{C^0} \\
\leq & C(f)\left\|g\right\|_{C^0}\leq C(f)\left\|g\right\|_{C^1},
\end{split}
\end{equation}
where $C(f)$ depends on the $C^0$ norm of $f^{-1}, \ f^{-2}$, and $f^{-3}$. 

When $\left\|f-f_0\right\|_{C^1}\leq \alpha$, we can show
\begin{equation}
\left\|f^{-1}\right\|_{C^0} \leq \left\|\left(f_0^{-1}(f-f_0)+I\right)^{-1}\right\|_{C^0}\left\|f_0^{-1}\right\|_{C^0} \leq C(\left\|f_0^{-1}\right\|_{C^0})
\end{equation}
as long as $\left\|f_0^{-1}\right\|_{C^0}\alpha< 1$. We can show similar results for $\left\|f^{-2}\right\|_{C^0}$ and $\left\|f^{-3}\right\|_{C^0}$. Therefore, for $f\in Z$,
\begin{equation}
\left\|DF_f(g)\right\|_{C^0} \leq C(f_0)\left\|g\right\|_{C^1}.
\end{equation}

Since
\begin{equation}
F(f)-F(f_0) = \int_0^1 \frac{d}{dt}F(f_0+t(f-f_0))\, dt = \int_0^1 DF_{f_0+t(f-f_0)}(f-f_0)\, dt, 
\end{equation}
we have
\begin{equation}
\left\|F(f)-F(f_0)\right\|_{C^0} \leq \int_0^1 \left\|DF_{f_0+t(f-f_0)}\right\|_{\mathcal{L}(C^1, C^0)}\left\|f-f_0\right\|_{C^1}\, dt \leq C(f_0)\left\|f-f_0\right\|_{C^1}, 
\end{equation}
which implies that $F:C^1(M)\rightarrow C^0(M)$ is locally Lipschitz. 

From \cite[Chapter 15, Proposition 1.1A]{taylor1997partial}, $f_t$ exists in short time.\\

\subsection{Induced Equations of Energy Density $e(f)$}

As for the Hessian term $\text{Hess}(W)_{AB}\langle df^A, df^B\rangle =\langle d\nabla W(f), df\rangle$, we have
\begin{equation}
\begin{split}
\langle d\nabla W(f), df\rangle = & \langle d(-f^{-3}), df\rangle \\
= & \langle f^{-3}df\, f^{-1}, df\rangle + \langle f^{-2}df\, f^{-2}, df\rangle + \langle f^{-1}df\, f^{-3}, df\rangle \\
= & |f^{-2}df|^2+2\langle f^{-2}df, f^{-1}df\, f^{-1}\rangle.
\end{split}
\end{equation}

When we study the induced Laplace or heat equation of the energy density, we have
\begin{equation}
\begin{split}
\Delta e(f) = & \Delta \left(\frac{1}{2}|df|^2+\frac{1}{2}|f^{-1}|^2\right)\\
= & -|\nabla df|^2-|f^{-3}|^2-\text{Ric}(df, df)-2|f^{-2}df|^2-4\langle f^{-2}df, f^{-1}df\, f^{-1}\rangle \\
\leq & |\text{Ric}|_\infty |df|^2+CW^2|df|^2\leq C_1e(f)+C_2e^3(f),
\end{split}
\end{equation}

and
\begin{equation}
\begin{split}
& \left(\frac{\partial }{\partial t}+\Delta\right)e(f_t)\\
 = & \left(\frac{\partial}{\partial t}+\Delta\right)\left(\frac{1}{2}|df_t|^2+\frac{1}{2}|f_t^{-1}|^2\right)\\
= & -|\nabla df_t|^2-|f_t^{-3}|^2-\text{Ric}(df_t, df_t)-2|f_t^{-2}df_t|^2-4\langle f_t^{-2}df_t, f_t^{-1}df_t\, f_t^{-1}\rangle\\
\leq & C_1e(f_t)+C_2e^3(f_t).
\end{split}
\end{equation}

Since the right-hand sides of the inequality in the last steps of both cases involve $e^3(f)$, we cannot simply apply the previous $\epsilon$-regularity to get the bound of the energy density $e$.

\section{Smoothing Potential $W_b$}

Since the unboundedness of $W$ gives undesired bounds for the Laplace or heat equation of the energy density $e$, we start with a bounded version of repulsive potential by defining
\begin{equation}
W_b(f) = \frac{1}{2}\text{trace}\left((f^2+bI)^{-1}\right)=\frac{1}{2}\sum_i^l\frac{1}{\lambda_i^2+b},
\end{equation}
where $I$ is the $l\times l$ identity matrix.

For $b>0$, we then analyze the bound of $e(f)$ in terms of $b$ using $\epsilon$-regularity. Later, we will show how the approximate co-dimension of the singularity set depends on $b$. 

Since 
\begin{equation}
W_b(f) \leq \frac{l}{2b}, \qquad |\text{Hess}(W_b)| \leq \frac{C(l)}{b}W_b(f) \leq \frac{C(l)}{b^2},
\end{equation}
we have the estimate of the Laplace and heat equation of the energy density $e_b$
\begin{equation}
\begin{split}
\Delta e_b(f) = &\Delta \left(\frac{1}{2}|df|^2+\frac{1}{2}\text{trace}\left((f^2+bI)^{-1}\right)\right)\\
 \leq &|\text{Ric}|_\infty|df|^2+\frac{C(l)}{b}W|df|^2\leq \left(\frac{C(l)}{b}e_b+C_2\right)e_b,
\end{split}
\end{equation}
and 
\begin{equation}
\Delta e_b(f) \leq \left(C_1+\frac{C_2(l)}{b^2}\right)e_b.
\end{equation}

Similarly,
\begin{equation}
\left(\frac{\partial}{\partial t}+\Delta \right)e_b(f_t) \leq \left(\frac{C(l)}{b}e+C_2\right)e_b,
\end{equation}
and
\begin{equation}
\left(\frac{\partial}{\partial t}+\Delta\right)e_b(f_t) \leq \left(C_1+\frac{C_2(l)}{b^2}\right)e_b.
\end{equation}

Here $C_1$ only depends on $M$, $C_2(l)$ depends on $l$, the size of $f$. 

Since $\Delta e_b$ or $\dsst \left(\frac{\partial }{\partial t}+\Delta\right)e_b$ is bounded by $e_b$, we immediately get the following result on the regularity of $e_b$:
\begin{theorem}[Elliptic]
	Fix a constant $b>0$. Let $f:M\rightarrow \R^{l\times l}\supset \text{Gr}(k, \mathds{C}^l)$ satisfy the following equation
	\begin{equation}\label{ellipticWb}
	d^*df+W_b(f) = 0,
	\end{equation}
	with $\dsst E_b(f)=\int_M \frac{1}{2}|df|^2+W_b(f)\, d\text{v}$ bounded. Then $e_b(f)=\dsst \frac{1}{2}|df|^2+W_b(f)\in C^\infty(M)$.
\end{theorem}
\begin{theorem}[Parabolic]
	Fix a constant $b>0$. Let $f_t:M\times I \rightarrow \R^{l\times l}\supset \textrm{Gr}(k, \mathds{C}^l)$ satisfy the equation
	\begin{equation}\label{paraWb}
	\frac{\partial f_t}{\partial t}+d^*df_t+W_b(f_t) = 0,
	\end{equation}
	with $f(0, x)= f_0(x)\in \text{Gr}(k, \mathds{C}^l)$ and smooth. Then $e_b(f_t)\in C^\infty(M)$.
\end{theorem} 

\begin{proof}
	Here we only prove the elliptic case. Since $E_b(f) = \dsst \int_M e_b(f)\, d\text{vol}$ is bounded, $e_b(f)\in L^1(M)$. That is, $|df|^2, \ W_b(f) \in L^1(M)$. Recall the inequality that
	$$\Delta e_b(f)\leq \frac{C}{b^2}e_b(f). $$
	
	From Moser iteration, we have
	\begin{equation}
	e_b(f)(x) \leq \left(\frac{C}{b^2}+\frac{C_2}{R^2}\right)^{m/2}\int_{B_R(x)} e(f)\, d\text{vol},
	\end{equation}
	for any $x\in M$. 
	
	Since $b$ is fixed, we can take $R\rightarrow\infty$ and get
	\begin{equation}
	e_b(f)(x) \leq C(b)\int_M e_b(f)\, d\textrm{vol} = C(b)E_b(f),
	\end{equation}
	which is bounded. Therefore, $e_b(f) \in L^\infty$ and $\Delta e_b(f) \in L^\infty$. Here we get $e_b(f)\in C^{1, \alpha}(M)$ for some $\alpha>0$. $e_b(f)\in C^\infty$ can be easily derived by bootstraping. 
\end{proof}

Since $e_b(f)\in C^\infty$, we consider $\sup e_b(f)$ and can get the following bound:
\begin{equation}
\sup_M e_b(f) \leq \left(\frac{C(l)}{b^2}\right)^{m/2}\int_M e_b(f)\, d\text{vol} \leq Cb^{-m}\int_M e_b(f)\, d\text{vol} = Cb^{-m}E_b(f).
\end{equation}

\subsection{Improved Bound of $\dsst \sup_M e_b$}

In the previous section, we have obtained two bounds for $-\textrm{Hess}(W_b)_{AB}$. They are either
\begin{equation}
-\textrm{Hess}(W_b)_{AB}\leq c(l)b^{-2}e_b(f),
\end{equation}
or
\begin{equation}
-\textrm{Hess}(W_b)_{AB}\leq c(l)b^{-1}\left(e_b(f)\right)^2.
\end{equation}

We can get an intermediate bound of $-\textrm{Hess}(W_b)_{AB}$ between these two
\begin{equation}
-\textrm{Hess}(W_b)_{AB} \leq c(l)b^{-1-q}\left(e_b(f)\right)^{2-q},
\end{equation}
for $0\leq q\leq 1$.

Using this intermediate bound, we can get finer control of $\dsst \sup_M e_b(f)$. In the following, we denote $\dsst \sup_M e_b(f)=e_M$. We illustrate the idea using the elliptic version.

We have in $M$,
\begin{equation}\label{intermediateinequa}
\begin{split}
\Delta e_b(f) = &-|\nabla df|^2-|\nabla W_b|^2-\textrm{Ric}(df, df)-2\textrm{Hess}(W_b)_{AB}\langle df^A, df^B\rangle \\
\leq & c(l)b^{-1-q}\left(e_b(f)\right)^{2-q}.
\end{split}
\end{equation}

We have shown that for a fixed $b>0$, $e_b(f)\in C^\infty(M)$, and therefore $e_M$ exists. From \eqref{intermediateinequa}, we get
\begin{equation}
\Delta e_b(f)\leq \left(c(l)b^{-1-q}e_M^{1-q}\right)e_b(f).
\end{equation}

Moser iteration gives
\begin{equation}\label{inter1}
e_M \leq \left(\frac{C_1(m)c(l)e_M^{1-q}}{b^{1+q}}+\frac{C_2(m)}{R^2}\right)^{\frac{m}{2}}\int_{B_R(x_0)} e_b(f)\, d\textrm{vol}.
\end{equation}
Here $e(x_0)=e_M$ and $R\leq 1$. In the following, for notational simplicity, we will assume $R_M\geq 1$. 

Multiplying $R^{\frac{2}{1-q}}b^{-\frac{1+q}{1-q}}$ on both sides and factoring out $1-q$ inside the parentheses on the right-hand side of \eqref{inter1}, we get
\begin{equation}
\begin{split}
\frac{R^{\frac{2}{1-q}}e_M}{b^{\frac{1+q}{1-q}}}\leq & \left(C_1\left(\frac{R^{\frac{2}{1-q}}e_M}{b^{\frac{1+q}{1-q}}}\right)^{1-q}+C_2\right)^{\frac{m}{2}}\left(\frac{R^{\frac{2}{1-q}}}{b^{\frac{1+q}{1-q}}}\right)\frac{1}{R^m}\int_{B_R(x_0)} e_b(f)\, d\textrm{vol}\\
= & \left(C_1\left(\frac{R^{\frac{2}{1-q}}e_M}{b^{\frac{1+q}{1-q}}}\right)^{1-q}+C_2\right)^{\frac{m}{2}}\frac{R^{p_0+\frac{2q}{1-q}}}{b^{\frac{1+q}{1-q}}}\frac{1}{R^{m-2+p_0}}\int_{B_R(x_0)} e_b(f)\, d\textrm{vol}\\
\leq & \left(C_1\left(\frac{R^{\frac{2}{1-q}}e_M}{b^{\frac{1+q}{1-q}}}\right)^{1-q}+C_2\right)^{\frac{m}{2}}\frac{R^{p_0+\frac{2q}{1-q}}}{b^{\frac{1+q}{1-q}}}\int_{B_1(x_0)} e_b(f)\, d\textrm{vol}\\
\leq & \left(C_1\left(\frac{R^{\frac{2}{1-q}}e_M}{b^{\frac{1+q}{1-q}}}\right)^{1-q}+C_2\right)^{\frac{m}{2}}\frac{R^{p_0+\frac{2q}{1-q}}}{b^{\frac{1+q}{1-q}}}\int_M e_b(f)\, d\textrm{vol}\\
= & \left(C_1\left(\frac{R^{\frac{2}{1-q}}e_M}{b^{\frac{1+q}{1-q}}}\right)^{1-q}+C_2\right)^{\frac{m}{2}}\frac{R^{p_0+\frac{2q}{1-q}}}{b^{\frac{1+q}{1-q}}} E_b(f).
\end{split}
\end{equation}
Here, we use the monotonicity formula and the fact that $R\leq 1$. 

Suppose we can pick $R$ so that 
\begin{equation}\label{star}
\dsst \frac{R^{\frac{2}{1-q}}e_M}{b^{\frac{1+q}{1-q}}}=1.
\end{equation}

This means
\begin{equation}
\begin{split}
R = & \left(e_Mb^{\frac{1+q}{1-q}}\right)^{\frac{1-q}{2}}=e_M^{\frac{1-q}{2}}b^{\frac{1+q}{2}}\\
\leq & \left(CE_b(f)b^{-m}\right)^{\frac{1-q}{2}}b^{\frac{1+q}{2}}\\
= & CE_b(f)^{-m(1-q)/2}b^{\frac{m+1}{2}q-\frac{m-1}{2}}.
\end{split}
\end{equation}

When $\dsst q>\frac{m-1}{m+1}$, we can guarantee the power $\dsst \frac{m+1}{2}q-\frac{m-1}{2}>0$, which ensures $R$ is small when $b<1$. 

Thus we pick $\dsst q=\frac{m-1}{m}=1-\frac{1}{m}>\frac{m-1}{m+1}$. We have
\begin{equation}
1 \leq (C_1+C_2)^{\frac{m}{2}}\frac{R^{p_0+2(m-1)}}{b^{2m-1}} E_b(f)=C(m)\frac{R^{p_0+2(m-1)}}{b^{2m-1}} E_b(f).
\end{equation}

This gives 
\begin{equation}
R\geq \left(\frac{b^{2m-1}}{C(m)E_b(f)}\right)^{\frac{1}{p_0+2(m-1)}}. 
\end{equation}

Recall we pick $R$ so that \eqref{star} holds. We get
\begin{equation}
e_M = \frac{b^{\frac{1+q}{1-q}}}{R^{\frac{2}{1-q}}} = \frac{b^{2m-1}}{R^{2m}}\leq b^{2m-1}\left(\frac{C(m)E_b(f)}{b^{2m-1}}\right)^{\frac{2m}{p_0+2(m-1)}}\leq C(m)E_b(f)^{\frac{2m}{p_0+2m-2}}b^{\frac{(2m-1)(p_0-2)}{p_0+2m-2}}.
\end{equation}
 
In the special case when $p_0=0$, we have
\begin{equation}\label{eMbound}
\sup_M e_b(f) = e_M \leq C(m)E_b(f)^{\frac{2m}{2m-2}}b^{-\frac{2m-1}{m-1}},
\end{equation}
which gives a better bound than $e_M\leq C(m)b^{-m}E_b(f)$. 

If \eqref{star} does not hold for any $R\leq R_M$, we have
\begin{equation}
\sup_M e_b(f)= e_M\leq \frac{b^{2m-1}}{R_M^{2m}}\leq b^{2m-1}.
\end{equation}

\subsection{$\epsilon$-Regularity Results}

Now we apply the $\epsilon$-regularity and get the bound of $e(f)$ in both cases in terms of the fixed $b$:

\begin{theorem}[Elliptic]
	Suppose $f:M\rightarrow R^{l\times l}$ is a smooth solution to \eqref{ellipticWb}. There exists $\epsilon_0 = \frac{b}{2C}$ for $C$ depend only on $M$ and $l$, such that if there exists some $x_0\in M$ and $R<R_M$ with $B_R(x_0)\subset M$ such that,  
	\begin{equation}
	\Phi(x_0, R) = \frac{1}{R^{m-2+p_0}}\int_{B_R(x_0)} (1+\nu(r, x))e_b(f)\, d\textup{vol} \leq \epsilon_0,
	\end{equation}
	then we have the estimate
	\begin{equation}
	\sup_{B_{\delta R}(x_0)}e(f)=\sup_{B_{\delta R}(x_0)}\left(\frac{1}{2}|df|^2+W_b(f)\right)\leq \frac{c(l)b}{(\delta R)^{2-p_0}},
	\end{equation}
	with $p_0$ as defined before, and $\delta \leq 3/4$.
\end{theorem}

\begin{proof}
	The energy density $\dsst e(f)=\frac{1}{2}|df|^2+W_b(f)$ satisfies the inequality in the ball $B_{\rho_0}(x_1)$, where $\rho_0$, $x_1$ are defined as in the previous $\epsilon$-regularity proof. 
	\begin{equation}
	\Delta e_b(f) \leq \left(\frac{C(l)}{b}\left(\sup_{B_{\rho_0}(x_1)}e_b\right)+C_2\right)e_b,
	\end{equation} 
	with 
	$$ \sup_{B_{\rho_0}(x_1)} e_b(f)\leq 4\sup_{B_{\sigma_0}(x_0)}e_b(f). $$
	
	Let 
	\begin{equation}
	e_0 = \sup_{B_{\delta R}(x_0)}e_b(f)=\sup_{B_{\delta R}(x_0)}\left(\frac{1}{2}|df|^2+W_b(f)\right) = e_b(x_1).
	\end{equation}
	
	Then we have
	\begin{equation}
	\Delta e(f)\leq \left(\frac{C(l)}{b}e_0+C_2\right)e\leq \frac{\tilde{C}(l)}{b}e_0 e(f).
	\end{equation}
	
	We can apply Lemma \ref{monoelliptic} and Theorem  \ref{epsilonelliptic} to get 
	\begin{equation}
	\rho^{2-p_0}e_0 \leq \left(\frac{C(l)}{b}\rho^{2-p_0}e_0+C_2\right)\frac{1}{\rho_0^{m-2+p_0}}\int_{B_{\rho_0}(x_1)}e_b(f)\, d\text{vol} \leq \left(\frac{C(l)}{b}\rho^{2-p_0}e_0+C_2\right)\epsilon_0.
	\end{equation}
	
	Suppose $C(l)\rho_0^{2-p_0}e_0/b = 1$. Picking $\epsilon_0 = \dsst \frac{b}{2C(M, l)}$ where
	$$ C(M, l) = C(l)(1+C_2(M))^{m/2} $$
	can lead to a contradiction. Therefore,
	\begin{equation}
	\sup_{B_{\delta R}(x_0)} e_b(f) \leq \frac{b}{C(l)(\delta R)^{2-p_0}} = \frac{c(l)b}{(\delta R)^{2-p_0}},
	\end{equation}
	with $C=\dsst \frac{\tilde{C}(l)}{b}$ and $c(l)=\dsst \frac{1}{\tilde{C}(l)}$. 
\end{proof}

We have the similar result for the parabolic case
\begin{theorem}[Parabolic]
	Suppose $f_t:M\times [-T, 0]\rightarrow \R^{l\times l}$ is a regular solution to parabolic equation \eqref{paraWb}. 
	
	There exists $\epsilon_0 = \dsst \frac{b}{2C}$, where $C$ only depends on $M$, and $l$, such that if for some $z_0=(x_0, t_0)\in M \times [-T,0]$, $R<R_M$, we have
	\begin{equation}
	\Psi(z_0, R) = \int_{t_0-4R^2}^{t_0-R^2}\frac{1}{(t_0-t)^{\nu_0}}\int_M G_{z_0}(x, t)e_b(f)\varphi^2\, d\textup{vol}\, dt\leq \epsilon_0,  
	\end{equation}
	with $P_R(z_0)\in M\times [-T,0]$, then we have the estimate
	\begin{equation}
	\sup_{P_{\delta R}(z_0)} e(f) = \sup_{P_{\delta R}(z_0)}\left(\frac{1}{2}|df_t|^2+W_b(f_t)\right) \leq \frac{c(l)b}{\left(\delta R\right)^{2-2\nu_0}},
	\end{equation}
	with $\delta$ only depending on $M$, $\inf\{R, 1\}$. 
\end{theorem}

Proof is similar to the elliptic case. 

\subsection{Measure of Set with Large $e_b(f)$}

\begin{definition}
	We define the Hausdorff $d$-measure at scale $b$ to be
\begin{equation}\label{hausdorff}
H_b^d(S) = \inf\left\{ \sum_{i=1}^{\infty} r_i^d\, \Big| \bigcup_{i=1}^\infty B_{r_i}(x_i)\supset S, r_i\leq b \right\}.
\end{equation}
\end{definition}

In this part we discuss the estimate of the Hausdorff measure at scale $b$ of the set where $e_b(f)$ is large. In the previous section, we have seen that $e_b(f)$ is $C^\infty$ with $b>0$ fixed. In order to analyze the effect when $b\rightarrow0$, we are interested in the size of the set
\begin{equation}\label{sigmabdef}
\Sigma_b=\left\{x\in M| e(f)(x)\geq \frac{1}{b} \right\}.
\end{equation}

Since $e_b(f)$ is at least continuous and $\Sigma_b$ is a closed set in a compact manifold $M$, $\Sigma_b$ is compact.

Here we only discuss the case when $p_0(x)=0$ and state the elliptic case. We have similar results with the parabolic metrics.

\begin{theorem}[Elliptic]\label{ellipticbadsetmeas}
	Fix a constant $b>0$. Let $f$ be the solution to \eqref{ellipticWb} with finite energy $E(f)$. Define the set $\Sigma_b$ as in \eqref{sigmabdef}. The Hausdorff $(m{-}1)$-measure of $\Sigma_b$ at scale $b$ is bounded and is independent of $b$
	\begin{equation}
	H^{m-1}_b(\Sigma_b) \leq CE(f).
	\end{equation}
\end{theorem}

\begin{proof}
From $\epsilon$-regularity, there exists $\epsilon_0=\dsst \frac{b}{2C}$, for $C$ only depends on $m$ and $M$, such that, if there exists $R$, 
\begin{equation}
\frac{1}{R^{m-2}}\int_{B_R(x_0)} (1+\nu(r,x))e(f)\, d\textrm{vol}\leq \epsilon_0=\frac{b}{2C},
\end{equation}
then \begin{equation}
\sup_{B_{R/2}(x_0)} e(f) \leq \frac{4c(l)b}{R^{2}}.
\end{equation}
Here $c(l)$ is an absolute constant depending on $l$, the size of $f$.

For $x_0\in \Sigma_b$, this implies
\begin{equation}
\frac{1}{R^{m-2}}\int_{B_R(x_0)}(1+\nu(r,x))e(f)\, d\textrm{vol}>\frac{b}{2C}\qquad \forall R>\sqrt{4c(l)}b.
\end{equation}

Let $\dsst \left\{B_{r_j}(x_j) \right\}_j$ be a cover of $\Sigma_b$ with $r_j=2\sqrt{(4c(l)}b$. Also assume $b$ is small enough so that $r_j<R_M/3$. By Vitali's covering lemma, we have finite sub-cover, still denoted as $\dsst \left\{B_{r_j}(x_j) \right\}_{j=1}^J$ such that $x_j\in \Sigma_b$, $B_{r_j}(x_j)$ are disjoint and $\dsst \bigcup_{j=1}^J B_{3r_j}(x_j)\supset \Sigma_b$. Inside each ball, we have
\begin{equation}
\frac{1}{r_j^{m-2}}\int_{B_{r_j}(x_j)} e(f)\, d\textrm{vol} > \frac{b}{2C}.
\end{equation}

That is
\begin{equation}
2C\int_{B_{r_j}(x_j)} (1+v(r, x))e(f)\, d\textrm{vol} > br_j^{m-2} .
\end{equation}

Summing over all $j$, we get
\begin{equation}
\begin{split}
3^{m-1}\sum_j^J br_j^{m-2}= & \sum_j^J (3r_j)^{m-1} \\
< & 2C(3^{m-1})\int_{\bigcup B_{r_j}(x_j)} (1+\nu(r, x))e_b(f)\, d\textrm{vol} \\
\leq &  C\int_M (1+\nu(r, x))e_b(f)\, dx \leq 2C\int_M e_b(f)\, d\textrm{vol}\\ \leq  & CE(f) .
\end{split}
\end{equation}

Since we can pick the initial map $f_0$ such that $E(f_0)$ is bounded and $f$ minimizes the energy functional $E$, $E(f)$ is finite. We get
\begin{equation}
H_b^{m-1}(\Sigma_b) \leq \sum_j^J(3r_j)^{m-1}<CE(f),
\end{equation}
which is finite.
\end{proof}

Heuristically, Theorem \ref{ellipticbadsetmeas} implies that as $b\rightarrow 0$, the dimension of singularity set of $f$ will be at most $m{-}1$.

\section{Higher Power Potentials}

We can see from $W_b(f)$ that the smaller $b$ is, the larger $e_b(f)$ can be. However, since $\epsilon_0$ depends on $b$,  $\epsilon_0\rightarrow0$ when $b\rightarrow 0$. The $\epsilon$-regularity results cannot be directly applied when $b\rightarrow0$. We consider changing powers of $f$ in the potential and define
\begin{equation}
W_b^L(f) = \frac{1}{2L}\text{trace}\left((f^{2L}+bI)^{-1}\right),\,\qquad L\geq 1.
\end{equation}
We hope to get tighter estimates for the expected co-dimension of the limiting singular set. By changing powers in the potential, we can get arbitrarily close to 2, suggesting that the earlier co-dimension 1 result is an artifact of an interaction of the specific potential and our techniques.

We consider the elliptic and parabolic equation for $f$ with this potential $W_b^L$
\begin{equation}\label{ellipticWbL}
d^*df+W_b^L(f) = 0 , \qquad \textrm{with $E_b^L(f)$ finite},
\end{equation}
and
\begin{equation}\label{paraWbL}
\frac{\partial f_t}{\partial t}+d^*df+W_b^L(f)=0,\ \qquad \textrm{with $f_0\in \textrm{Gr}(k, \mathds{C}^l)$ and $E(f_0)$ finite}.
\end{equation}

It is obvious that
\begin{equation}
W_b^L(f)\leq \frac{c(l)}{Lb}.
\end{equation}

Detailed computations give
\begin{equation}
\left(W_b^L\right)_B = -\textrm{trace}\left((f^{2L}+bI)^{-2}f^{2L-1}B\right),\qquad \nabla W_b^L = - \left(f^{2L}+bI\right)^{-2}f^{2L-1} .
\end{equation}

We write $f$ in a local orthonormal basis $\{e_i \}_{i=1}^l$
\begin{equation}
f = \sum_i^l\lambda_i e_i\otimes e_i^*. 
\end{equation}

Then
\begin{equation}
|df|^2= \sum_i|d\lambda_i|^2 + \sum_{i,j}(\lambda_i-\lambda_j)^2|\langle de_i, e_j\rangle |^2,
\end{equation}
and
\begin{equation}
\nabla W_b^L(f) = -\sum_i (\lambda_i^{2L}+b)^{-2}\lambda_i^{2L-1}e_i\otimes e_i^*.
\end{equation}

The Hessian term of $\langle d\nabla W_b^L, df\rangle$ can be computed as follows
\begin{equation}\label{higherHessian}
\begin{split}
  &   \langle d\nabla W_b^L, df\rangle \\
= & -\sum_i\langle d(\frac{\lambda_i^{2L-1}}{(\lambda_i^{2L}+b)^{2}})  , d\lambda_i \rangle\\
& -\sum_{i,j}\langle (\frac{\lambda_i^{2L-1}}{(\lambda_i^{2L}+b)^{2}}-\frac{\lambda_j^{2L-1}}{(\lambda_j^{2L}+b)^{2}})   (\lambda_i-\lambda_j)|\langle d e_i,e_j\rangle |^2\\
= & -\sum_i (2L-1) \frac{\lambda_i^{2L-2}}{(\lambda_i^{2L}+b)^{2}}|d\lambda_i|^2+4 L  \frac{\lambda_i^{4L-2}}{(\lambda_i^{2L}+b)^{3}} |d\lambda_i |^2\\
  & -\sum_{i,j}(\lambda_i-\lambda_j)^{-1} \left(\frac{\lambda_i^{2L-1}}{(\lambda_i^{2L}+b)^{2}}-\frac{\lambda_j^{2L-1}}{(\lambda_j^{2L}+b)^{2}}\right)   (\lambda_i-\lambda_j)^2|\langle d e_i,e_j\rangle |^2\\
= &
\sum_i\lambda_i^{2L-2}\frac{(2L+1)\lambda_i^{2L}-(2L-1)b}{(\lambda_i^{2L}+b)^{3}} |d\lambda_i |^2\\
 & -\sum_{i,j}(\lambda_i-\lambda_j)^{-1}  \left(\frac{\lambda_i^{2L-1}}{(\lambda_i^{2L}+b)^{2}}-\frac{\lambda_j^{2L-1}}{(\lambda_j^{2L}+b)^{2}}\right)(\lambda_i-\lambda_j)^2|\langle d e_i,e_j\rangle |^2.
\end{split}
\end{equation}

Observe that when $|\lambda_i|>b^{\frac{1}{2L}}$, the first term $\dsst \lambda_i^{2L-2}\frac{(2L+1)\lambda_i^{2L}-(2L-1)b}{(\lambda_i^{2L}+b)^{3}} |d\lambda_i |^2$ is positive. When $|\lambda_i<b^{\frac{1}{2L}}|$, we have
\begin{equation}
\left|\lambda_i^{2L-2}\frac{(2L+1)\lambda_i^{2L}-(2L-1)b}{(\lambda_i^{2L}+b)^{3}}\right| \leq c(L)b^{1-\frac{1}{L}}b^{-2}.
\end{equation}

For the second term in \eqref{higherHessian}, denote
$$ g(\lambda) = \frac{\lambda^{2L-1}}{(\lambda^{2L}+b)^2} .$$
Then
\begin{equation}
\begin{split}
 & (\lambda_i-\lambda_j)^{-1}  \left(\frac{\lambda_i^{2L-1}}{(\lambda_i^{2L}+b)^{2}}-\frac{\lambda_j^{2L-1}}{(\lambda_j^{2L}+b)^{2}}\right)\\
= & \frac{g(\lambda_i)-g(\lambda_j)}{\lambda_i-\lambda_j}\\
\leq & C|g'(\lambda)|, \qquad \textrm{for some $\lambda$ between $\lambda_i$ and $\lambda_j$}.
\end{split}
\end{equation}

Since 
$$ g'(\lambda) = \lambda^{2L-2}\frac{(2L+1)\lambda^{2L}-(2L-1)b}{(\lambda^{2L}+b)^3}, $$
we have shown previously that when $|\lambda|<b^{\frac{1}{2L}}$, $$|g'(\lambda)|\leq c(L)b^{-1-\frac{1}{L}}$$

When $|\lambda|\geq b^{1/2L}$, we have
$$ |g'(\lambda)| \leq c(L)\lambda^{-2L-2}\leq c(L)b^{-1-\frac{1}{L}}. $$

Therefore, in general we always have
\begin{equation}
\left|\langle d\nabla W_b^L, df\rangle\right| \leq c(l, L)b^{-1-\frac{1}{L}}|df|^2.
\end{equation}

Since $W_b^L\leq \dsst \frac{c(l, L)}{b}$, we can conclude that $e_b^L(f)\in C^\infty(M)$. We have two versions to control the term $-2\langle d\nabla W_b^L, df\rangle$ inside the induced elliptic or parabolic equation of $e_b^L(f)$, which will give us different controls of $e_b^L(f)$. We illustrate the idea using the elliptic equation.

When using
\begin{equation}
\begin{split}
\Delta e_b^L(f)= &-|\nabla df|^2-|\nabla W_b^L|^2-\textrm{Ric}(df, df)-2\langle d\nabla W_b^L, df\rangle\\
        \leq & c(l,L)b^{-1-\frac{1}{L}}e_b^L,
\end{split}
\end{equation}
we get the upper bound of $\sup e$ by using Moser iteration directly
\begin{equation}
\sup_M e(f) \leq c(l, L)b^{-(1+\frac{1}{L})\frac{m}{2}}\int_M e(f)\, d\textrm{vol} = c(l, L)b^{-(1+\frac{1}{L})\frac{m}{2}} E(f).
\end{equation}

\subsection{$\epsilon$-Regularity Results}

When we want to estimate the Hausdorff measure at scale $b^{\frac{1}{2}+\frac{1}{2L}}$ of the set
\begin{equation}
\Sigma_b^L = \left\{x\in M | e_b^L(f)\geq \frac{1}{b} \right\},
\end{equation}
we need a finer estimate of $e_b^L$
\begin{equation}
\begin{split}
\Delta e_b^L(f) = & -|\nabla df|^2-|\nabla W_b^L|^2-\textrm{Ric}(df, df)-2\langle d\nabla W_b^L, df\rangle\\
\leq & c(l, L)b^{-\frac{1}{L}}(e_b^L)^2.
\end{split}
\end{equation}

Since the monotonicity formula (Lemma \ref{monoelliptic}) still holds for $e_b^L$, we derive the $\epsilon$-regularity results

\begin{theorem}[Elliptic]
	Suppose $f$ solves the elliptic equation \eqref{ellipticWbL} with finite energy $E_b^L(f)$. 
	
	There exists $\epsilon_0 = \dsst\frac{b^{\frac{1}{L}}}{2C}$ with $C$ depending only on $M, l, L$, such that if there exists $R<\rho$, and
	\begin{equation}
	\Phi_b^L(x_0, R) = \frac{1}{R^{m-2+p_0}}\int_{B_R(x_0)} (1+\nu(r, x))e_b^L(f)\, d\textup{vol} \leq \epsilon_0,
	\end{equation}
	then for any $\delta \leq \frac{3}{4}$, 
	\begin{equation}
	\sup_{B_{\delta R}(x_0)} e_b^L(f) \leq \frac{c(L, l, M)b^{\frac{1}{L}}}{(\delta R)^{2-p_0}}.
	\end{equation}
	Here, $\nu(r, x)$ and $p_0$ are defined as before.
\end{theorem}

\begin{proof}
	The key to the proof lies in the elliptic inequality
	\begin{equation}
	\Delta e_b^L(f) \leq c(l, L)b^{-\frac{1}{L}}(e_b^L)^2.
	\end{equation}
	
	We apply Moser iteration to get
	\begin{equation}
	\sup_{B_{\delta R}(x_0)} e_b^L \leq \left(\frac{C_1(M, l, L)}{b^{\frac{1}{L}}}+\frac{C_2(M)}{(R-\delta R)^2}\right)^{\frac{m}{2}}\int_{B_R(x_0)} e_b^L(f)\, d\textrm{vol}.
	\end{equation}
	
	This implies that the function $h(\sigma)$ satisfies
	\begin{equation}
	h(\sigma_0) \leq \left(\frac{C_1(M, l, L)h(\sigma_0)}{b^{\frac{1}{L}}}+C_2(M)\right)^{m/2}\frac{1}{R^{m-2+p_0}}\int_{B_R(x_0)} e_b^L(f)\, d\textrm{vol}.
	\end{equation}
	
	We can pick $\dsst \epsilon_0 \leq \frac{b^{\frac{1}{L}}}{2C(M, l, L)}$, where $C(M, l, L)$ is determined by $C_1(M, l, L)$ and $C_2(M)$. Then a contradiction argument implies 
	$$ h(\sigma) \leq h(\sigma_0) \leq 1.$$
	
	That is,
	\begin{equation}
	\sup_{B_{\delta R}(x_0)} e_b^L(f) \leq \frac{c(M, l, L)b^{\frac{1}{L}}}{(\delta R)^{2-p_0}}.
	\end{equation}
\end{proof}

\begin{theorem}[Parabolic]
	Suppose $f_t$ solves the parabolic equation \eqref{paraWbL} with finite initial energy $E_b^L(f_0)$.
	
	There exists $\epsilon_0=\dsst \frac{b^{\frac{1}{L}}}{2C}$, with $C$ depending only on $M, l, L$, such that if there exists $R$ with $0<R<\sqrt{t_0-t}/2$ and
	\begin{equation}
	\Psi((x_0, t_0), R) = \int_{t_0-4R^2}^{t_0-R^2}\frac{1}{(t_0-t)^{\nu_0}}\int_M G_{(x_0, t_0)}(x, t)e_b^L(f_t)\varphi^2\, d\textup{vol}\, dt \leq \epsilon_0,
	\end{equation}
	then for any $\delta$ depending only on $M, R$, we have
	\begin{equation}
	\sup_{P_{\delta R}(z_0)} e_b^L(f) \leq \frac{c(L, l, M)b^{\frac{1}{L}}}{(\delta R)^{2-2\nu_0}}.
	\end{equation} 
	Here $z_0=(x_0,t_0)$, $\nu_0$ is the same ratio defined in \eqref{mu0nu0} and the theorem holds when $\nu_0=0$. 
\end{theorem}

We omit the proof here since it is an easy analog of the elliptic case. 

\subsection{Measure of Set with Large $e_b^L(f)$}

In this part we assume $p_0=0$. 

We consider the set
\begin{equation}\label{sigmabLdef}
\Sigma_b^L = \left\{x\in M | e_b^L(f)\geq \frac{1}{b} \right\},
\end{equation}

Since $e_b^L(f)$ is continuous and $\Sigma_b^L$ is a closed set in a compact manifold $M$, $\Sigma_b^L$ is compact.

\begin{theorem}[Elliptic]
	Fix a constant $b>0$. Let $f$ be the solution to \eqref{ellipticWbL}. Define the set $\Sigma_b^L$ as in \eqref{sigmabLdef}. The Hausdorff $\left(\frac{2}{L+1}+(m-2)\right)$-measure of $\Sigma_b^L$ at scale $b^{\frac{1}{2}+\frac{1}{2L}}$ is bounded and independent of $b$.
\end{theorem}

\begin{proof}
	From $\epsilon$-regularity, there exists $\epsilon_0=\dsst \frac{b^{\frac{1}{L}}}{2C}$, for $C$ only depending on $M, l, L$, such that, if there exists $R$, 
	\begin{equation}
	\frac{1}{R^{m-2}}\int_{B_R(x_0)} (1+\nu(r,x))e_b^L(f)\, d\textrm{vol}\leq \epsilon_0=\frac{b^{\frac{1}{L}}}{2C},
	\end{equation}
	then 
	\begin{equation}
	\sup_{B_{R/2}(x_0)} e_b^L(f) \leq \frac{4c(M, l, L)b^{\frac{1}{L}}}{R^{2}}.
	\end{equation}
	
	For $x_0\in \Sigma_b$, this implies
	\begin{equation}
	\frac{1}{R^{m-2}}\int_{B_R(x_0)}(1+\nu(r,x))e(f)\, d\textup{vol}>\frac{b^{\frac{1}{L}}}{2C}\qquad \forall R>\sqrt{4cb^{1+\frac{1}{L}}}
	\end{equation}
	
	Let $\dsst \left\{B_{r_j}(x_j) \right\}_j$ be a cover of $\Sigma_b^L$ with $r_j=2\sqrt{(4cb^{1+\frac{1}{L}}}$. By Vitali's covering lemma, we have a finite sub-cover, still denoted as $\dsst \left\{B_{r_j}(x_j) \right\}_{j=1}^J$ such that $x_j\in \Sigma_b^L$, $B_{r_j}(x_j)$ are disjoint and $\dsst \bigcup_{j=1}^J B_{3r_j}(x_j)\supset \Sigma_b^L$. Inside each ball, we have
	\begin{equation}
	\frac{1}{r_j^{m-2}}\int_{B_{r_j}(x_j)} e_b^L(f)\, d\textrm{vol} > \frac{b^{\frac{1}{L}}}{2C}.
	\end{equation}
	That is
	\begin{equation}
	2C\int_{B_{r_j}(x_j)} (1+v(r, x))e_b^L(f)\, d\textrm{vol} > 2C\int_{B_{r_j}(x_j)}e_b^L(f)\, d\textrm{vol} > b^{\frac{1}{L}}r_j^{m-2} .
	\end{equation}
	Summing over all $j$, we get
	\begin{equation}
	\begin{split}
	3^{\frac{2}{L+1}+(m-2)}\sum_j^J b^{\frac{1}{L}}r_j^{m-2}= & \sum_j^J (3r_j)^{\frac{2}{L+1}+(m-2)} \\
	< & 2C(3^{\frac{2}{L+1}+(m-2)})\int_{\bigcup B_{r_j}(x_j)} (1+\nu(r, x))e_b^L(f)\, d\textrm{vol} \\
	\leq &  C\int_M (1+\nu(r, x))e_b^L(f)\, dx \leq 2C\int_M e(f)\, d\textrm{vol}\\ \leq  & CE_b^L(f_0) .
	\end{split}
	\end{equation}
	
	Since we can pick $f_0$ such that $E(f_0)$ is bounded and $f$ minimizes the energy functional $E$, we get
	\begin{equation}
	H_b^{(\frac{2}{L+1}+(m-2)}(\Sigma_b^L) \leq \sum_j^J(3r_j)^{\frac{2}{L+1}+(m-2)}<CE_b^L(f_0).
	\end{equation}
	which is finite. Therefore, the Hausdorff $\frac{2}{L+1}+(m-2)$-measure of $\Sigma_b^L$ at scale $b^{\frac{1+L}{2L}}$ is bounded. 
\end{proof}

